\documentclass[12pt]{article}
\usepackage{amsmath, amsfonts, amsthm, amssymb, verbatim}
\usepackage{graphicx}
\usepackage[mathcal]{euscript}
\usepackage{amstext}
\usepackage{amssymb,mathrsfs}
\usepackage{bbm}
\usepackage{color}
\usepackage[latin1]{inputenc}
\usepackage{txfonts}
\newcommand{\clg}[1]{{\mathcal{#1}}}

\newcommand{\R}{\mathbb R}

\newcommand{\C}{\mathbb C}

\newcommand{\ve}{\varepsilon}

\newcommand{\sgn}{\text{sgn}}

\newcommand \loc    {\text{loc}}

\newtheorem{theorem}{Theorem}[section]
\newtheorem{proposition}[theorem]{Proposition}

\newtheorem{lemma}[theorem]{Lemma}

\newtheorem{definition}[theorem]{Definition}

\begin{document}
\title{On Fractional Benney Type Systems}

\author{Wladimir Neves$^1$, Dionicio Orlando$^1$}

\date{}

\maketitle

\footnotetext[1]{Mathematical Institute, Universidade Federal
do Rio de Janeiro, C.P. 68530, Cidade Universit\'aria 21945-970,
Rio de Janeiro, Brazil. E-mail: {\sl wladimir@im.ufrj.br,
domorenov@unac.edu.pe.}}


%
%
\begin{abstract} This paper introduces fractional type evolutionary equations modeling the interaction between 
short waves and long waves. We consider a fractional Benney type system, which is given by a fractional Schr\"odinger equation
coupled with a fractional porous medium equation. 
Under the assumption of weak coupling or small initial data related to the fractional Schr\"odinger equation, it is proved
the existence of weak solutions to the Cauchy problem. 
\end{abstract}

\maketitle

\section {Introduction} \label{Intro}

The main issue of this paper is to introduce 
and study the Cauchy problem for the following fractional Benney type systems
\begin{equation}
\label{ecua00}
\left\{
\begin{aligned}
 i \ \partial_t u - (- \Delta)^{s} u & = \alpha \ v \ u + \gamma \ |u|^{2}u,  \hspace{10pt} \ x\in \R, \ t>0,
  \\[5pt]
\partial_t v + (-\Delta)^{s/2} g(v) & = \beta \ (-\Delta)^{s/2} |u|^{2},  \hspace{18pt} \ x\in \R, \ t> 0,
  \\[5pt]
u(0,x)=u_{0}(x), \ \ &v(0,x)=v_{0}(x),  \hspace{22pt} \ x \in \R, 
\end{aligned}
\right.
\end{equation}
where $\alpha$, $\beta$ and $\gamma$ are real constants. 
The complex value function $u(t,x)$ is the unknown of the fractional Schr\"odinger equation, 
which describes the short wave, and the real value function $v(t,x)$ is the unknown of the fractional 
porous medium equation, which describes the long wave.
Here $(- \Delta)^{s}$, $(0 < s < 1)$, denotes the usual fractional Laplacian in $\R^n$, which characterize 
nonlocal, long-range diffusion effects and can be defined by 
$\mathcal{F}\{(- \Delta)^{s}f\}(\xi) = |\xi|^{2s}\mathcal{F}\{f\}(\xi)$, where $\mathcal{F}$ is the Fourier Transform. 
The function $g \in C^{1}(\R)$ is assumed to be nondecreasing,
hence degenerated zones for the state variable $v(t,x)$ are allowed. 
A particular case of \eqref{ecua00}, (e.g. $g \equiv$ constant),
it has its own interest
 \begin{equation}
\label{ecua01}
\left\{
\begin{aligned}
 i \ \partial_t u - (- \Delta)^{s} u & = \alpha \ v \ u +\gamma \ |u|^{2}u,  \hspace{7pt} \ x\in \R, \ t> 0,
  \\[5pt]
\partial_t v & = \beta \ (-\Delta)^{s/2} |u|^{2},  \hspace{18pt} x \in \R, \ t> 0,
  \\[5pt]
u(0,x)=u_{0}(x), \ \ &v(0,x)=v_{0}(x),  \hspace{23pt} x \in \R.
\end{aligned}
\right.
\end{equation}

The theory of evolutionary equations modeling the interaction between 
short waves and long waves goes back to Benney \cite{B:Benney}.
Indeed, in that paper Benney propose a general system (see equations (3.27), (3.28) in that paper), and 
we recall below the closer one studied by Bekiranov, Ogawa, Ponce \cite{BOP}, that is to say 

\begin{equation}\label{ecua02}
\left\{
\begin{aligned}
 &i \ \partial_t S - (-\Delta) S + i \ C_{S } \ \nabla S  = \alpha SL + \gamma \ |S|^{2}S,   \hspace{11pt} x \in \R, \ \ t> 0, 
 \\[5pt] 
& \partial_t L + C_{L} \ \nabla L + \nu P(D_x)L + \lambda \ \nabla L^2 = \beta \ \nabla |S|^{2},  \hspace{8pt} x\in \R, \ \ t> 0,
\end{aligned}
\right.
\end{equation}
where $C_{S,} \alpha, \gamma, C_{L}, \nu, \lambda$ and $\beta$ are real constants. 
Moreover, $P(D_x)$ is a linear differential operator with constant coefficients.
Applying a proper gauge transformation and a scaling of the variables, 
the system \eqref{ecua02}, when $\nu=0$, is equivalent to

\begin{equation}
\label{ecua03}
\left\{
\begin{aligned}
 i \ \partial_t u - (-\Delta) u & = \alpha \ v \ u + \gamma \ |u|^{2}u, 
  \\[5pt]
\partial_t v + C \ \nabla v & = \beta \ \nabla |u|^{2},
\end{aligned}
\right.
\end{equation}
where $C = \pm 1$. In fact, the authors in \cite{BOP} claim that, the system \eqref{ecua03}
is the most typical case in the theory of wave interaction.

In particular, for $s= 1$ and $g(v)= v$, the system \eqref{ecua00} recalls \eqref{ecua03}, since we have the following equivalence
$$
   \| (-\Delta)^{1/2}f \|_{L^{2}(\R^n)}= \| \nabla f\|_{L^{2}(\R^n)}, \quad \text{for each $f \in H^1(\R^n)$}.
$$ 
Then, one may roughly speaking interpret 
 \eqref{ecua00} as a generalization of \eqref{ecua03}.
 In particular, the system \eqref{ecua00} makes sense for $x \in \R^n$ and $t> 0$.  
Although, this is not exactly the case. 
Indeed, even if $\nabla f$ and $ (-\Delta)^{1/2} f$ have the same $L^2-$norm
they are different objects, that is, the former has local behavior and the other is
nonlocal.  

Therefore, we highlight the motivations to
consider the fractional Benney type systems proposed in this paper, besides the multidimensional 
one. Indeed, the short 
(transversal) wave described by the Schr\"odinger equation may represent a signal (wave packets), that is
$u(t,x)$ is a function that conveys information to control, for instance, 
some underwater equipment.
This information propagates in a generalized medium, where long (longitudinal) waves 
are described by the porous medium equation. Here,
the fractional Laplacian introduces the long-range interactions in 
both equations, which are coupled by the $\alpha$, $\beta$ constants.
%
This discussion follows to applications in 
Synthetic Aperture Radar (see \cite{Beal}), and 
atmospheric internal gravity waves (see \cite{Tabaei}, \cite{Wilhelm}), 
which represent complex anomalous
systems and it seems better modeled by fractional Laplacians. 
%
%

Last but not least, Benney in \cite{B:Benney} also consider the following system (see equations $(3.8)$, $(3.9)$ in that paper)
\begin{equation}
\label{ecua04}
\left\{
\begin{aligned}
iS_{t} - (-\Delta) S + i \ C_{g} \ \nabla S &= \alpha \ LS+ \gamma \ |S|^{2} S,  \hspace{9pt} \ x\in \R, \ \ t> 0, 
\\[5pt]
L_{t} + C_{l} \ \nabla L&= \beta \ \nabla |S|^{2},  \hspace{45pt} \ x\in \R, \ \ t> 0.
\end{aligned}
\right.
\end{equation}
In particular, when $C_{g}= C_{l}$ long waves and short waves are resonant, and in this case Tsutsumi and Hatano in \cite{TsutHat1} 
proved that, the transformation: $x \mapsto y= x - C_{g} \ t$ eliminates the first $x$-derivative terms in \eqref{ecua04}, 
hence we have
\begin{equation}
\label{ecua05}
\left\{
\begin{aligned}
iu_{t} - (-\Delta) u &= \alpha \ v u+ \gamma \ |u|^{2}u,
\\[5pt]
v_{t}&= \beta \ \nabla|u|^{2},
\end{aligned}
\right.
\end{equation}
which resembles the fractional short wave and long wave system 
\eqref{ecua01}. 

\bigskip
{\bf \large Statement of the Main Result.} 

The following definition tells us in which sense a pair $(u(t,x),v(t,x))$ is a 
weak solution to the Cauchy problem \eqref{ecua00}.
Hereafter, we fix $\gamma= 1$, and without loss of generality $ g(0)= 0$.
\begin{definition}
\label{DEFSOL}
Given an initial data $(u_{0}, v_{0}) \in H^{s}(\R) \times (L^{2}(\R) \cap L^{\infty}(\R))$
and any $T> 0$, a pair 
$
 (u,v) \in L^{\infty}(0,T; H^{s}(\R)) \times L^{\infty}(0,T; L^{2}(\R) \cap L^{\infty}(\R))
$
is called a weak solution of the Cauchy problem \eqref{ecua00},
when it satisfies:
\begin{equation}
\label{RP1}
\begin{aligned}
    & i \!\! \int_{0}^{T} \int_{\R} \Big( u(t,x) \ \partial_t \overline{\varphi}(t,x) 
    +  (-\Delta)^{s/2} u(t,x) \ (-\Delta)^{s/2} \overline{\varphi}(t,x) \Big)  dx dt 
    + i \!\! \int_{\R} u_{0}(x) \ \overline{\varphi}(0,x) dx
\\[5pt]
    &\quad + \alpha  \int_{0}^{T}\int_{\R} v(t,x) \ u(t,x) \ \overline{\varphi}(t,x) dx dt 
    +  \int_{0}^{T} \int_{\R} |u(t,x)|^{2} \ u(t,x) \ \overline{\varphi}(t,x) dx dt= 0,
\end{aligned}
\end{equation}
\begin{equation}
\label{RP2}
\begin{aligned}
   \int_{0}^{T} \int_{\R} v(t,x) \ \partial_t \psi(t,x)
   &- g(v(t,x)) \ (-\Delta)^{s/2} \psi(t,x) dx dt 
   + \int_{\R} v_{0}(x) \ \psi(0,x) dx
   \\[5pt]
    &+ \beta \int_{0}^{T} \int_{\R} |u|^{2}(t,x) \ (-\Delta)^{s/2} \psi(t,x) dxdt= 0, 
\end{aligned}
\end{equation}
for each test function $\varphi, \psi \in C^{\infty}_{c} \big( (-\infty,T) \times \R \big)$, 
with $\varphi$ being complex-valued and $\psi$ real-valued.
\end{definition}

Now, we state plainly the main
result of this paper.
\begin{theorem}[Main Theorem]
\label{MAINTHM}
Let $(u_{0}, v_{0}) \in H^{s}(\R) \times (L^{2}(\R) \cap L^{\infty}(\R))$, $(\frac{1}{2} < s < 1)$,
and $g \in C^1(\R)$ satisfying 
$$
    0 \leq g'(\cdot) \leq M < \infty. 
$$
For any $T> 0$, there exist $\alpha_0> 0$, $E_0> 0$, such that, if 
$|\alpha| \leq \alpha_0$ or $ \|u_{0}\|_{L^{2}(\R)} \leq E_0$, then
there exists a weak solution 
$$
 (u,v) \in L^{\infty}(0,T; H^{s}(\R)) \times L^{\infty}(0,T; L^{2}(\R) \cap L^{\infty}(\R))
$$
of the Cauchy problem \eqref{ecua00}. Moreover, for a.a. $t \in (0, T)$
\begin{equation}
\label{maxprinciple}
   \|v(t)\|_{L^\infty(\R)} \leq \|v_0\|_{L^\infty(\R)}. 
\end{equation}

\end{theorem}
Clearly, how lower 
are the $\alpha, \beta$ constants less coupled are the equations in \eqref{ecua00}. 
In fact, the $\alpha$ constant makes the
difference concerning the global in time existence, (see Theorem \ref{MAINTHM}). 
Another very important point is the energy input to the
signal, i.e. $\|u_0\|_{L^2}$. As far as the information has to be sent, more energy 
is needed. Again, the statement of the Main Theorem shows that, the global in time 
solvability depends on the amount of energy given to the signal. 
 
Finally, we recall that the fractional Schr\"odinger equation
appears in the water wave models in \cite{IONESCU}. In fact,
the fractional Schr\"odinger equation was introduced in the theory 
related to fractional quantum mechanics associated to $s$-stable 
L\'evy process (see for instance \cite{Laskin}). This field is developing fast, hence jointly with 
\cite{IONESCU} we address the 
reader to the following papers \cite{ChoHajaiejHwangOzawa}, 
\cite{ChoHwangKwonLee} and \cite{GuoHuo}. 
Moreover, the fractional
porous medium equations has been widely studied 
in the last years. For instance, we address V\'azquez \cite{JLV} (and references there in),
where is described the physical and mathematical background related to 
nonlinear diffusion equations involving nonlocal effects. 

\section {Notation and Background} \label{preliminary}

In this section we fix the notations, and collect some preliminary results. 
First, let $\Omega \subset \R^n$ be open set. 
We denote by $dx, d\xi$, etc. 
the Lebesgue measure on $\Omega$ and by $L^p(\Omega)$, $p \in [1,+\infty)$, the set of (real or complex) $p$-summable functions 
with respect to the Lebesgue measure.
%
Moreover, we denote by
$\mathcal{F} \varphi(\xi) \equiv \widehat{\varphi}(\xi)$
the Fourier Transform of $\varphi$, which is an
isometry in $L^2(\R^n)$.


$\bullet$ {\bf The space $W^{s,p}(\Omega)$}

The Sobolev space is denoted by $W^{s,p}(\Omega)$, where a real $p \geqslant 1$
is the integrability index and a real $s \geqslant 0$ is the
smoothness index. More precisely, for $s \in (0,1)$, 
$p \in [1,+\infty)$, the fractional Sobolev space of
order $s$ with Lebesgue exponent $p$ is defined by
$$
W^{s,p}(\Omega) := \Big\{ u \in L^p(\Omega) : \int_{\Omega} \int_{\Omega} \dfrac{\vert u(x)-u(y)\vert^p}{
\vert x-y\vert ^{n+sp}} \ dx \ dy < + \infty \Big\},
$$
endowed with norm
$$
\Vert u\Vert_{W^{s,p}(\Omega)} = \left(  \int_{\Omega} \vert u\vert^p dx + \int_{\Omega} \int_{\Omega} \dfrac{\vert u(x)-u(y)\vert^p}{
\vert x-y\vert ^{n+sp}} \ dx \ dy \right)^\frac{1}{p}.
$$
For $s > 1$ we write $s = m + \sigma$, where $m$ is an integer
and $\sigma \in (0, 1$). In this case, the space $W^{s,p}(\Omega)$ consists of those equivalence classes
of functions $u \in W^{m,p}(\Omega)$ whose distributional derivatives $D^{\alpha} u$, with $|\alpha| = m$,
belong to $W^{\sigma,p}(\Omega)$, that is
$$
W^{s,p}(\Omega) = \Big\{ u \in W^{m,p}(\Omega) : 
{\sum_{\vert \alpha\vert = m}}\Vert D^{\alpha}u \Vert_{W^{\sigma,p}(\Omega)} < \infty \Big\},
$$
which is a Banach space with respect to the norm
$$
\Vert u\Vert_{W^{s,p}(\Omega)} = \Big(  \Vert u\Vert^p_{W^{m,p}(\Omega)} 
+ {\sum_{\vert u\vert = m}}\Vert D^{\alpha}u \Vert^p_{W^{\sigma,p}(\Omega)} \Big)^\frac{1}{p}.
$$
If $s= m$ is an integer, then the space $W^{s,p}(\Omega)$ coincides with the Sobolev space
$W^{m,p}(\Omega)$. It is very interesting the case when $p = 2$, i.e. $W^{s,2}(\Omega)$,
which is also a Hilbert space and we can consider the inner product
$$
\langle u, v\rangle_{W^{s,2}(\Omega)} = \langle u,v\rangle + \int_{\Omega} \int_{\Omega} 
\frac{(u(x)-u(y))}{\vert x-y\vert ^{\frac{n}{2}+s}} \ \frac{(v(x)-v(y))}{\vert x-y\vert ^{\frac{n}{2}+s}} \ dx \ dy,
$$
where $\langle\cdot,\cdot\rangle$ is the inner product in $L^2(\Omega)$.

\bigskip
$\bullet$ {\bf The space $H^s(\R^n)$}

Now, following Tartar \cite{Tartar} 
we take into account an alternative definition of the space $H^s(\R^n)= W^{s,2}(\R^n)$ via
Fourier Transform. Precisely, we may define

\begin{equation}
\label{HsDEF}
    H^s(\R^n):= \Big\{ u \in L^2(\R^n) : \int_{\R^{n}} (1+|\xi|^{2s}) \ |\mathcal{F} u(\xi)|^{2} \ d\xi < \infty \Big\}
\end{equation}
and we observe that the above definition, is valid also for any real $s \geq1$. Moreover, 
$H^s(\R^n)$ is a Hilbert space with the scalar product
$$
   (u,v)_{H^{s}(\R^{n})} = \int_{\R^{n}} (1+|\xi|^{2})^{s} \ \widehat{u}(\xi) \ \overline{\widehat{v}}(\xi) \ d\xi.
$$

The equivalence of the above definitions is stated in the following
\begin{lemma}
Let $0 < s < 1$. Then, the definitions of $H^s(\R^n)$ and $W^{s,2}(\R^n)$ are equivalent. In particular, for any 
$u \in H^s(\R^n)$
 
\begin{equation}
\label{SOVFRACFT}
       \int_{\R^n} \int_{\R^n} \frac{|u(x)-u(y)|^{2}}{|x-y|^{n+2s}} \ dx \ dy = 2 \ C_{n,s}^{-1} \ 
       \int_{\R^n} |\xi|^{2 s} \ |\mathcal{F}u(\xi)|^{2} d\xi,
\end{equation}
where 
$$
    C_{n,s}^{-1} = \int_{\R^{n}} \frac{1-\cos(\zeta_{1})}{|\zeta|^{n+2s}} \ d\zeta.
$$
\end{lemma}

One remarks that, for $s> n/2$, the Hilbert space $H^s(\R^n)$ is an algebra (see \cite{FP:Felipe}). Moreover, there exists a constant $C= C(s)> 0$, such that 
for any $f, g \in  H^s(\R^n)$
\begin{equation}
\label{algebra}
    \|f \ g\|_{H^s(\R^n)} \leq C  \ \|f\|_{H^s(\R^n)}  \|g\|_{H^s(\R^n)}.
\end{equation}

\subsection{Fractional Laplacian operator in $\R^n$}

The fractional Laplacian operator can be defined in $\R^n$ by
\begin{equation}
\label{LAPFRACFT1}
 \widehat{(-\Delta)^{s} \ f}(\xi) = |\xi|^{2s} \ \hat{f}(\xi), \qquad (0<s<1). 
\end{equation}
Hence the fractional Laplacian is a pseudo-differential
operator with principal symbol $|\xi|^{2s}$.
The fractional Laplacian 
can be similarly described using singular integrals 
\begin{equation}
\label{LAPFRACFT2}
 (-\Delta)^s f(x)= C_{n,s} \, {\rm P.V.}\! \int_{\mathbb
R^n} \frac{f(x)-f(\xi)}{|x-\xi|^{n+2s}} \ d\xi.
\end{equation}
Moreover, its inverse denoted by $\clg{K}_s:= (-\Delta)^{-s}$, $(0 < s < 1)$,
is given by convolution with the Riesz kernel $K_s(x) = C_{n,s} \ |x|^{2s-n}$, that is,
$\clg{K}_s f= K_s \ast f$. 

It follows from \eqref {HsDEF}, \eqref{SOVFRACFT} and \eqref{LAPFRACFT1} that, 
there exist positive constants $m_s$, $M_s$, such that, for each $f \in H^s(\R^n)$
\begin{equation}
\label{NORMHs}
m_s \big(\|f\|_{L^2(\R^n)} + \| (-\Delta)^{s/2} f\|_{L^2(\R^n)}\big) \leq
\|f\|_{H^s(\R^n)} \leq M_s \big(\|f\|_{L^2(\R^n)} + \| (-\Delta)^{s/2} f\|_{L^2(\R^n)}\big).
\end{equation}

$\bullet$ {\bf Bilinear form} 

In order to study the fractional
diffusion term, it will be important to associate a bilinear form
to the operator $\mathcal{K}_s$ in the space $H^s(\mathbb R^n),$
$0 < s < 1,$ which is given for any pair $v,w \in H^s(\mathbb R^n)$ by
\begin{equation}
\label{BWV}
\mathcal B_s(v,w): = C_{n,s} \iint_{\mathbb
R^{2n}} \left( v(x)-v(y) \right) \frac{1}{|x-y|^{n+2s}} \left( w(x)-w(y) \right) \ dx dy.
\end{equation}
The bilinear form
$\mathcal B_s$ were considered in \cite{caff-vazquez-soria-jems} as
an auxiliary tool in the study of regularity properties of
solutions to the fractional type porous medium equation. 

\begin{lemma}[See \cite{caff-vazquez-soria-jems}]
\label{lem.bilform}
If $v$ is given by $v = G(w),$ with $G' \geq 0,$ then,
$\mathcal B_s(v,w) \geq 0$. 
Furthermore, for every $v,w \in H^1(\mathbb R^n)$ we have the
characterization
\begin{equation}
\label{BWV10}
\mathcal B_s(v,w)= C \iint_{\mathbb R^{2n}} \nabla
v(x) \frac{1}{|x-y|^{n-2+2s}} \nabla w(y) \ dx dy,
\end{equation}
where $C$ is a positive constant. 
\end{lemma}

\begin{proposition}\label{BWV11}
Let $v \in H^{1}(\R^{n})$, $G \in C^1(\R)$ with $G'(\cdot) \geq m > 0$. Then
$$
\int_{\R^n} (-\Delta)^{s/2} G(v) \ v \ dx \geq m \ C_{n,s}^{-1} \ \|(-\Delta)^{s/4}v\|_{L^{2}(\R^n)}^{2}.
$$
\end{proposition}
\begin{proof}
It follows directly from \eqref{BWV}, \eqref{BWV10} and applying the
intermediate value theorem. 
\end{proof}
\subsection{Auxiliary kernels}

$\bullet$ {\bf Unitary group for the Schr\"odinger equation}

For each $\ve > 0$, we consider the following Cauchy problem for $u(t,x) \in \C$, driven by the linear 
fractional perturbed Schr\"odinger equation
\begin{equation}
\label{LSEQ}
\left\{
\begin{aligned}
 i \ \partial_{t} u-(-\Delta)^{s}u -\ve^{a} (-\Delta)u = 0, &  ~~ x \in \R^n ~~ t \in \R, 
\\[5pt]
 u(0,x) = u_{0}(x), &  ~~ x\in \R^n,
\end{aligned}
\right.
\end{equation}
where $a \in \R$ is a fixed parameter chosen a posteriori. 
Applying the Fourier transform in the spatial variable, we have
\begin{equation*}
\left\{
\begin{aligned}
 i \ \partial_{t} \widehat{u}(t,\xi) - |\xi|^{2s} \ \widehat{u}(t,\xi) - \ve^{a} \ |\xi|^{2} \ \widehat{u}(t,\xi) = 0, &  ~~ \xi \in  \R^n ~~ t \in \R, 
 \\[5pt]
\widehat{u}(0,\xi) = \widehat{u}_{0}(\xi), &  ~~ \xi \in \R^n,
\end{aligned}
\right.
\end{equation*}
which solution is given by $\widehat{u}(t,\xi) = e^{- i \ \big( |\xi|^{2s} + \ve^{a} \ |\xi|^{2} \big) t} \widehat{u_{0}}(\xi)$.
Therefore, it follows that
$$
    u(t,x) = \mathcal{F}^{-1} \bigg\{ e^{- i \ \big( |\xi|^{2s} + \ve^{a} \ |\xi|^{2} \big) t}\mathcal{F}u_{0}(\xi) \bigg\}(x)
$$
solves the Cauchy problem \eqref{LSEQ}. 
For $u_{0} \in L^{2}(\R^n)$, ($\mathcal{F} u_{0} \in L^{2}(\R^n)$), then
$$
     e^{- i \ \big( |\xi|^{2s} + \ve^{a} \ |\xi|^{2} \big) t}\mathcal{F}u_{0}(\xi)\in L^{2}(\R^n).
$$
Now, we define for each $t \in \R$ the operator
\begin{equation}
\label{GROUPS}
u \mapsto  U_{\ve}(t) u:= \mathcal{F}^{-1} e^{- i \ \big( |\xi|^{2s} + \ve^{a} \ |\xi|^{2} \big) t}\mathcal{F}u,
\end{equation}
which is bounded in $L^2(\R^n)$ for each $u \in L^{2}(\R^n)$. Indeed, we have
$$
\begin{aligned}
  \|U_{\ve}(t) u\|_{L^{2}(\R^n)}^{2} & = \int_{\R^n} |U_{\ve}(t)u(x)|^{2} \ dx
   = \int_{\R^n} |\widehat{U_{\ve}(t)u}(\xi)|^{2} \ d\xi
    \\[5pt]
   & = \int_{\R^n} |e^{- i \ \big( |\xi|^{2s} + \ve^{a} \ |\xi|^{2} \big) t} \widehat{u}(\xi)|^{2} \ d\xi
    = \int_{\R^n} |\widehat{u}(\xi)|^{2} \ d\xi.
\end{aligned}
$$
Therefore, the family $(U_{\ve}(t))_{t \in \R}$ is a group of isometries in $L^2(\R^n)$. 

One remarks that, $H^{s}(\R^n)$, $(s> 0)$, is invariant by the isometry group $(U_{\ve}(t))_{t \in \R}$.
For each $u \in H^{s}(\R^n)$, we have
$$
\begin{aligned}
  \|U_{\ve}(t) u\|_{H^{s}(\R^n)} & = \int_{\R^n} (1 + |\xi|^{2})^{s} \ |\widehat{U_{\ve}(t) u}(\xi)|^{2} \ d\xi
    \\[5pt]
   & = \int_{\R^n} (1 + |\xi|^{2})^{s} |\widehat{u}(\xi)|^{2} \ d\xi
   =\|u\|_{H^{s}(\R^n)}.
\end{aligned}
$$
Thus $U_{\ve}(t)(H^{s}(\R^n))$ is a closed subspace in $H^{s}(\R^n)$ and, we have 
$$
    H^{s}(\R^n)= U_{\ve}(t)(H^{s}(\R^n)) \oplus (U_{\ve}(t)(H^{s}(\R^n)))^\perp.
$$    
Moreover, since $U_\ve(t)$ is symmetric in $H^s(\R^n)$
$$
\begin{aligned}
  (U_{\ve}(t)u,w)_{H^{s}(\R^n)} & = \int_{\R^n} (1 + |\xi|^{2})^{s} \widehat{U_{\ve}(t) u}(\xi) \ \widehat{w}(\xi) \ d\xi
   \\[5pt]
   & = \int_{\R^n} (1 + |\xi|^{2})^{s} e^{- i \ \big( |\xi|^{2s} + \ve^{a} \ |\xi|^{2} \big) t} \ \widehat{u}(\xi) \ \widehat{w}(\xi) \ d\xi
    \\[5pt]
   & =\int_{\R^n} (1 + |\xi|^{2})^{s} \ \widehat{u}(\xi) \ \widehat{U_{\ve}(t)w}(\xi) \ d\xi
    =(u,U_{\ve}(t) w)_{H^{s}(\R^n)}
\end{aligned}
$$
and also an isometry, it follows that $(U_{\ve}(t)(H^{s}(\R^n)))^\perp= \{0\}$. 

\bigskip
$\bullet$ {\bf Semigroups of contractions for the heat equation}

For each $\ve > 0$, we consider the following Cauchy problem for $v(t,x) \in \R$, driven by the linear 
Heat equation
\begin{equation}
\label{LHEQ}
 \left\{
  \begin{aligned}
   \partial_{t} v - \ve^{b} \Delta v = 0, &  ~~ x \in \R^n, ~~ t > 0, 
   \\[5pt]
    v(0,x) = v_{0}(x), &  ~~ x \in \R^n,
   \end{aligned}
\right.
\end{equation}
where $b \in \R$ is a fixed parameter chosen a posteriori. 
Again, applying the Fourier transform in the spatial variable, we obtain
\begin{equation*}
\left\{
\begin{aligned}
 \partial_{t} \widehat{v}(t,\xi) + \ve^{b} \ |\xi|^{2} \ \widehat{v}(t,\xi) = 0, &  ~~ \xi \in \R^n, ~~ t> 0, 
 \\[5pt]
 \widehat{v}(0,\xi) = \widehat{v}_{0}(\xi), &  ~~ \xi \in \R^n,
\end{aligned}
\right.
\end{equation*}
which solution is given by $\widehat{v}(t,\xi) = e^{- \ve^{b} \ |\xi|^{2} t} \ \widehat{v_{0}}(\xi)$. 
Consequently,
$$
    v(t,x) = \mathcal{F}^{-1} \bigg \{ e^{- \ \ve^{b} \ |\xi|^{2} t} \ \mathcal{F}v_{0}(\xi) \bigg\}(x)
$$
solves the Cauchy problem \eqref{LHEQ}, and it is well known that,
for $v_{0} \in L^{2}(\R^n)$, ($\mathcal{F}v_{0} \in L^{2}(\R^n)$), it follows that 
$e^{- \ve^{b} \ |\xi|^{2} t} \ \mathcal{F}v_{0}(\xi) \in  L^{2}(\R^n)$. 

Similarly, we define for each $t> 0$ the operator
\begin{equation}
\label{GROUPH}
  v \mapsto W_{\ve}(t) v = \mathcal{F}^{-1} e^{- \ve^{b} \ |\xi|^{2} t} \ \mathcal{F} v.
\end{equation}
The operator $W_{\ve}(t)$ is bounded in $L^2(\R^n)$, in fact 
the family $\{W_{\ve}(t) v\}_{t>0}$ is a semigroup of contractions.
Indeed, for any $t> 0$, 
$\|W_{\ve}(t) v_{0}\|_{L^{2}}
\leq \|v_{0}\|_{L^{2}}$,
for any $v \in L^{2}(\R^n)$.
Also in $H^{1}(\R^n)$, that is  
$$
\begin{aligned}
  \|W_{\ve}(t) v_{0}\|_{H^{1}}^{2}
   = \!\!\int_{\R^n} (1\!\! +\! |\xi|^{2}) \ | \widehat{W_{\ve}(t)v_{0}}(\xi) |^{2} \ d\xi
    = \!\!\int_{\R^n} (1\!\! +\! |\xi|^{2}) \ | e^{- \ve^{b} |\xi|^{2} t} \ \widehat{v_{0}}(\xi) |^{2} \ d\xi
   \leq  \|v_{0}\|^{2}_{H^{1}}.
\end{aligned}
$$

One recalls that, the Heat kernel has a regularity effect. Indeed, 
a refined estimate is given by the following 
\begin{lemma}\label{E02}
For any $v \in L^{2}(\R)$, there exists a constant $C>0$ independent of $t$ and $v$, such that for any  $t> 0$
\begin{equation}
\label{E03}
   \|\partial_{x} W_{\ve}(t)v\|_{L^{2}(\R)} \leq \frac{C}{\sqrt{t}} \|v\|_{L^{2}(\R)}.
\end{equation}
\end{lemma}

\begin{proof}
From \eqref{GROUPH}, we have
$$
\begin{aligned}
W_{\ve}(t)v &= \mathcal{F}^{-1} e^{- \ve^{b} \ |\xi|^{2} t} \ \mathcal{F} v
= \frac{1}{\sqrt{4\pi \ve^{b} t}} e^{- \frac{x^{2}}{4 \ve^{b}t}} \ast v,
\end{aligned}
$$
then
$$
\partial_{x} W_{\ve}(t) v= \frac{1}{\sqrt{4\pi \ve^{b} t}} \frac{-2x}{4\ve^{b}t} e^{- \frac{x^{2}}{4 \ve^{b}t}} \ast v.
$$
Applying Young's inequality, it follows that
$$
\begin{aligned}
   \|\partial_{x} W_{\ve}(t)v\|_{L^{2}(\R)} &\leq \frac{2}{4\ve^{b}t} \ 
    \|\frac{x}{\sqrt{4\pi \ve^{b} t}} \ e^{- \frac{x^{2}}{4 \ve^{b}t}}\|_{L^{1}(\R)} \|v\|_{L^{2}(\R)}
     =  \frac{1}{\sqrt{\pi \ve^{b}}} \ \frac{1}{t^{1/2}} \ \|v\|_{L^{2}(\R)}.   
\end{aligned}
$$
\end{proof}

\subsection{Auxiliary inequalities} 
\label{Ineq}

The next two auxiliary results will be used broadly in this paper.

\begin{proposition}
\label{Chain Rule}(Chain Rule) 
Let $f \in H^{s}(\R^n)$, $0 < s < 1$, $F \in C^{1}(\C)$ 
with $\|F^{\prime}\|_{L^{\infty}(\R)} \leq M$ for some $M > 0$. Then

\begin{equation}\label{ecua81}
\|(-\Delta)^{s/2} F(f)\|_{L^{2}(\R^n)} \leq \|F^{\prime}\|_{L^{\infty}(\R)} \ \|(-\Delta)^{s/2} f\|_{L^{2}(\R^n)}
\end{equation}
\end{proposition}

%

\medskip
Now, we provide the following (sharp) result. 
\begin{proposition}
\label{product} 
Let $f \in H^{s}(\R)$, $\frac{1}{2}< s < 1$. Then, 
\begin{equation}
\label{DESG1product}
  \|f\|_{L^{\infty}(\R)} \leq \frac{2}{\sqrt{\pi (2s-1)}} \ \|f\|_{L^{2}(\R)}^{1-\frac{1}{2s}} 
  \ \|(-\Delta)^{s/2} f\|_{L^{2}(\R)}^{\frac{1}{2s}}, 
\end{equation}
and
\begin{equation}
\label{DESG2product}
  \|(-\Delta)^{s/2} |f|^{2}\|_{L^{2}(\R)} \leq 2 
  \ \|f\|_{L^{\infty}(\R)} \ \|(-\Delta)^{s/2} f\|_{L^{2}(\R)}. 
\end{equation}
\end{proposition}

\begin{proof} 1. First, since $s > 1/2$, it follows from the well-known Embedding Theorem that, 
$H^s$ is an algebra of functions. Moreover, 
a function $f \in H^s(\R)$ may be represented by a continuous function which vanishes at infinity. 
Let us show \eqref{DESG1product}, hence applying the inverse Fourier transform, we have for each 
$x \in \R$ 

$$
\begin{aligned}
  |f(x)| & = \Big| \frac{1}{(2 \pi)^{1/2}} \int_{\R} e^{i x \xi}\widehat{f}(\xi) \ d\xi \Big| 
   \leq \frac{1}{(2 \pi)^{1/2}} \int_{\R} |\widehat{f}(\xi)| \ d\xi 
   \\[5pt]
   & = \frac{1}{(2 \pi)^{1/2}} \Big( \int_{|\xi|\leq R}| \widehat{f}(\xi)| \ d\xi
   + \int_{|\xi|\geq R} \frac{|\xi|^{s}}{|\xi|^{s}} |\widehat{f}(\xi)| \ d\xi \Big),
\end{aligned}
$$
where $R > 0$ is any fixed real number.
Then, applying the Cauchy-Schwartz inequality 

\begin{equation}
\label{ecua86}
\begin{aligned}
  |f(x)| & \leq \frac{1}{(2 \pi)^{1/2}}  \Big( \int_{|\xi| \leq R}1 \ d\xi \bigg)^{1/2} \bigg( \int_{|\xi|\leq R} |\widehat{f}(\xi)|^{2} \ d\xi \bigg)^{1/2} 
  \\[5pt]
   & + \frac{1}{(2 \pi)^{1/2}} \Big( \int_{|\xi|\geq R} \frac{1}{|\xi|^{2s}} \ d\xi \bigg)^{1/2}
   \bigg( \int_{|\xi|\geq R} |\xi|^{2 s} \ |\widehat{f}(\xi)|^{2} \ d\xi \bigg)^{1/2}
   \\[5pt]
   & \leq \frac{1}{(2 \pi)^{1/2}} \Big( \sqrt{2} \ R^{1/2} \ \|f\|_{L^{2}(\R)}
   + \sqrt{\frac{2}{2s-1}} \ R^{\frac{1}{2}-s} \ \|(-\Delta)^{s/2}f\|_{L^{2}(\R)} \Big)
   \\[5pt]
   & \leq \frac{1}{\sqrt{\pi}\sqrt{2s-1}} \Big( R^{1/2} \|f\|_{L^{2}(\R)} + R^{\frac{1}{2}-s} \|(-\Delta)^{s/2} f\|_{L^{2}(\R)}  \Big).
\end{aligned}
\end{equation}

Conveniently, we consider $R= \|f\|_{L^{2}(\R)}^{-\frac{1}{s}} \ \|(-\Delta)^{s/2} f\|_{L^{2}(\R)}^{\frac{1}{s}}$ in \eqref{ecua86} to obtain

$$
  |f(x)| \leq 
   \frac{1}{\sqrt{\pi (2s-1)}} \Big( \|f\|_{L^{2}(\R)}^{1-\frac{1}{2s}}\|(-\Delta)^{s/2} f\|_{L^{2}(\R)}^{\frac{1}{2s}} + \|f\|_{L^{2}(\R)}^{1-\frac{1}{2s}}\|(-\Delta)^{s/2} f\|_{L^{2}(\R)}^{\frac{1}{2s}} \Big).
$$

2. Now, we prove \eqref{DESG2product}. Again, from \eqref{LAPFRACFT2} and the definition of the Fractional Laplacian, we obtain
$$
\begin{aligned}
  &\|(-\Delta)^{s/2} |f|^{2}\|_{L^{2}(\R)}^{2}  = \frac{C_{n,s}}{2}\int_{\R} \int_{\R} \frac{||f|^{2}(x)-|f|^{2}(y)|^{2}}{|x-y|^{1+2s}} \ dx dy 
   \\[5pt]
   & \leq C_{n,s} \ \Big( \iint_{\R \times \R} \frac{|f(x) \ (\overline{f}(x) - \overline{f}(y))|^{2}}{|x-y|^{1+2s}} \ dx dy 
   + \iint_{\R \times \R}  \frac{|\overline{f}(y) \ (f(x) - f(y))|^{2}}{|x-y|^{1+2s}} \ dx dy \Big)
   \\[5pt]
   & \leq 2 \ \|f\|_{L^{\infty}(\R)} \ \|(-\Delta)^{s/2} f\|_{L^{2}(\R)}.
\end{aligned}
$$
\end{proof}

\subsection{Generalized Grownwall Lemma}
\label{GGL}
We consider the following (see \cite{SS:Sever}),
\begin{theorem}
Let $\eta(t)$ be a nonnegative function which satisfies the inequality

\begin{equation}\label{ecuaI18}
\eta(t) \leq C + \int_{t_{0}}^{t} \bigg( a(\tau) \ \eta(\tau) + b(\tau) \ \eta^{\sigma}(\tau) \bigg) \ d\tau, ~~ C \geq 0, \sigma \geq 0,
\end{equation}
where $a(t)$ and $b(t)$ are continuous nonnegative functions for $t \geq t_{0}$.

1. For $0 \leq \sigma < 1,$
\begin{equation}
\label{ecuaI19}
\begin{aligned}
     \eta(t) \leq \bigg\{ C^{1 - \sigma} \exp\bigg[ &(1 - \sigma)\int_{t_{0}}^{t} a(\tau) \ d\tau \bigg]
     \\
     + & (1 - \sigma) \int_{t_{0}}^{t} b(\tau) \exp \bigg[ (1 - \sigma) \int_{\tau}^{t} a(r) \ dr \bigg] \ d\tau \bigg\}^{\frac{1}{1-\sigma}}.
\end{aligned}
\end{equation}

2. For $\sigma = 1,$
\begin{equation}\label{ecuaI20}
\eta(t) \leq C \exp \bigg\{ \int_{t_{0}}^{t} \big[ a(\tau) + b(\tau) \big] \ d\tau \bigg\}.
\end{equation}

3. For $\sigma > 1$, with the additional hypothesis
\begin{equation}
\label{ecuaI21}
    C < \bigg\{ \exp \bigg[ (1 - \sigma) \int_{t_{0}}^{t_{0}+h} a(\tau) \ d\tau \bigg] \bigg\}^{\frac{1}{\sigma - 1}}
    \bigg\{ (\sigma - 1) \int_{t_{0}}^{t_{0}+h} b(\tau) \ d\tau \bigg\}^{- \frac{1}{\sigma - 1}},
\end{equation}
we also get for $t_{0} \leq t \leq t_{0}+h$, for $h > 0$
\begin{equation}\label{ecuaI22}
\eta(t) \leq C \ \Big\{ \exp \Big[ (1 - \sigma) \!\! 
\int_{t_{0}}^{t} a(\tau) \ d\tau \Big] - C^{-1} (\sigma - 1) \!\!
\int_{t_{0}}^{t} b(\tau) \exp \Big[ (1 - \sigma) \!\!
\int_{\tau}^{t} a(r) \ dr \Big] \ d\tau \Big\}^{\frac{1}{\sigma-1}}.
\end{equation}
\end{theorem}

\subsection{Entropies} 
\label{Entropies}

Following the scalar conservation laws theory, we say that a Lipschitz 
convex function $\eta: \R \to \R$ is an entropy. 
The most important example is the family of Kru\v zkov's entropies, that is 
$$
   \eta_k(v):= |v - k|, \quad \text{for each $k \in \R$}. 
$$
Then, we recall that any smooth entropy $\eta(v)$, which is linear at infinity, 
can be recovered by the family of  Kru\v zkov's entropies. Indeed, a straight calculation 
shows that 
$$
   \eta(v)= \frac{1}{2} \int_\R \eta^{\prime \prime}(\xi) \, |v - \xi| \, d\xi,
$$
modulo an additive constant. 
Symilarly, given $g \in C^1(\R)$ and 
$q: \R \to \R$, such that, $q^\prime= \eta^\prime \, g^\prime$, then 
$$
  q(v)= \frac{1}{2} \int_\R \eta^{\prime \prime}(\xi) \, |g(v) - g(\xi)| \, d\xi. 
$$ 
Under the above conditions, $(\eta,q)$ is called here an
entropy pair. 

\medskip
Now, we consider the following 
\begin{lemma}
Let $v$ be a real $H^{1}(\R)$ function, $g \in C^1(\R)$ satisfying 
$$
   0 \leq g^\prime(\cdot) \leq M < \infty,
$$
and $s \in (0,1)$. 
Then, for each $k \in \R$ fixed, and each $x \in \R$, 
\begin{equation}
\label{remainder} 
\begin{aligned}
(-\Delta)^{s}|g(v(x)) - g(k)|= \sgn(v(x)-k) \, (-\Delta)^{s} g(v(x)) - R_{k}(x), 
\end{aligned}
\end{equation} 
where the non-negative remainder function $R_k(\cdot)$ is given by 
\begin{equation}
\label{RK}
R_{k}(x):=\left\{
\begin{aligned}
    & 2 \, C_{1,s} \int_{\{v(y)<k\}} \frac{g(k) - g(v(y))}{|x-y|^{1+2s}} \ dy,  \quad \{v(x)>k\},
 \\[5pt]
   & 2 \, C_{1,s} \int_{\{v(y)>k\}} \frac{g(v(y))-g(k)}{|x-y|^{1+2s}} \ dy,  \quad \{v(x)<k\}.
\end{aligned}
\right.
\end{equation}
\end{lemma} 
\begin{proof} Since the function $g$ is non-decreasing, it follows that 
$$
   \sgn(v(x) - k) (g(v(x) - g(k))= | g(v(x) - g(k) |.
$$
Therefore, we have 
$$
\begin{aligned}
(-\Delta)^{s}&|g(v(x)) - g(k)| 
 = (-\Delta)^{s} \big((\sgn(v(x)-k)(g(v(x))-g(k))\big)
\\[5pt]
& = \sgn(v(x)-k)(-\Delta)^{s}g(v(x)) + (g(v(x))-g(k)) \, (-\Delta)^{s}\sgn(v(x)-k)
\\[5pt]
& - C_{1,s} \int_{\R}\frac{\big(\sgn(v(x)-k)-\sgn(v(y)-k)\big)}{|x-y|^{1+2s}} \big((g(v(x))-g(v(y)))\big) \ dy
\\[5pt]
& = \sgn(v(x)-k)(-\Delta)^{s} g(v(x))
\\[5pt]
& \ \ \ \ + (g(v(x)) -g(k)) C_{1,s} \int_{\R}\frac{\sgn(v(x)-k)-\sgn(v(y)-k)}{|x-y|^{1+2s}} \ dy
\\[5pt]
& \ \ \ \ - C_{1,s} \int_{\R} (g(v(x))-g(v(y)))\frac{\sgn(v(x)-k)-\sgn(v(y)-k)}{|x-y|^{1+2s}} \ dy
\\[5pt]
& = \sgn(v(x)-k)(-\Delta)^{s} g(v(x))
\\[5pt]
& \ \ \ \ + C_{1,s} \int_{\R} (g(v(y))-g(k))\frac{\sgn(v(x)-k)-\sgn(v(y)-k)}{|x-y|^{1+2s}} \ dy,
\end{aligned}
$$
where we have used that $\{x\in \R:v(x)>k\}$
and $\{x\in \R:v(x)<k\}$ are open sets, since $v$ is continuous. 
\end{proof}

\section {On a Perturbed System} 
\label{Localexistence}

In order to show 
the solvability of the Cauchy problem \eqref{ecua00},
we perturbe both equations $(\ref{ecua00})_1$, and 
$(\ref{ecua00})_2$, adding 
Laplacian terms with different 
velocities of perturbation.  
Specifically, let $a, b > 0$ be fixed parameters and for each $\ve \in (0,1)$, we consider the following
system posed in $(0,T) \times \R$, 
\begin{equation}
\label{ecuaI2}
\left\{
\begin{aligned}
    & i \ \partial_t u^{\ve} - (- \Delta)^{s} u^{\ve} + \ve^{a}  \Delta u^{\ve}= \alpha \ v^{\ve} \ u^{\ve} + |u^{\ve}|^{2}u^{\ve},  
 \\[5pt]
   & \partial_t v^{\ve} - \ve^{b} \ \Delta v^{\ve}= \beta \ (-\Delta)^{s/2} (|u^{\ve}|^{2}) - (- \Delta)^{s/2} g_\ve(v^{\ve}),  
 \\[5pt]
  &  u^{\ve}(0,x) = u^\ve_{0}(x), ~~ v^{\ve}(0,x) = v^\ve_{0}(x),
\end{aligned}
\right.
\end{equation}
where $T> 0$ is a real number, conveniently $g_\ve(v):= g(v) + \ve v$, and 
the pair $(u^\ve_0, v^{\ve}_0) \in H^{1}(\R)  \times H^{1}(\R)$ is an approaching sequence 
converging strongly to
$(u_0, v_0)$ in $H^{s}(\R)  \times (L^{2}(\R) \cap L^{\infty}(\R))$, ($\|v^\ve_0\|_{L^\infty(\R)} \leq \|v_0\|_{L^\infty(\R)}$). 
First, we show (local in time) existence  
and uniqueness of mild solution to \eqref{ecuaI2}.
Then, we derive a priori 
important estimates, which enable us to extend the local in time solution.
Moreover, we stress that these a priori estimates will be also important to 
show that the family $\{(u^\ve, v^\ve)\}$ of solution to \eqref{ecuaI2} is 
relatively compact. 

\subsection{Existence and uniqueness} 

The following definition tell us in which sense the pair $(u^{\ve},v^{\ve})$ is a solution
of the Cauchy problem \eqref{ecuaI2}.
\begin{definition}
\label{mildsolution}
The pair $(u^{\ve},v^{\ve}) \in C \big( [0,T];H^{1}(\R) \big) \times C \big( [0,T];H^{1}(\R) \big)$ 
is called a mild solution of \eqref{ecuaI2} if satisfies the following integral equations

\begin{equation}
\label{ecua152}
\left\{
\begin{aligned}
 & u^{\ve}(t) = U_{\ve}(t) \ u^\ve_{0} - i \int_{0}^{t} U_{\ve}(t-t^{\prime}) \bigg( \alpha \ v^{\ve}(t^{\prime}) \ u^{\ve}(t^{\prime})
 + \ |u^{\ve}(t^{\prime})|^{2} \ u^{\ve}(t^{\prime}) \bigg) \ dt^{\prime}, 
   \\[5pt]
& v^{\ve}(t) = W_{\ve}(t) \ v^\ve_{0} + \int_{0}^{t} W_{\ve} \big( t - t^{\prime} \big) \bigg( \beta \ (-\Delta)^{s/2} |u^{\ve}(t^{\prime})|^{2}
- (-\Delta)^{s/2} g_\ve(v^{\ve}) \bigg) \ dt^{\prime},
\end{aligned}
\right.
\end{equation}
where $U_{\ve}(t)$, $W_{\ve}(t)$ are given respectively by 
\eqref{GROUPS} and \eqref{GROUPH}.  
\end{definition}

We are going to apply the Banach Fixed Point Theorem to show the local-in-time existence of solutions as defined above. 
To begin, we consider the following lemma (we put $\ve= 1$ for simplicity with obvious notation).

\begin{lemma}\label{15L1}
Let $\frac{1}{2} < s < 1$, $g \in C^{1}(\R)$, satisfying $1 \leq g^{\prime}(\cdot) \leq M < \infty$, $(g(0)=0)$. 
For $T > 0$, let $(\tilde{u},\tilde{v}) \in C \big( [0,T];H^{1}(\R) \big) \times C \big( [0,T];H^{1}(\R) \big)$, then for each
$(u_{0},v_{0}) \in H^{1}(\R) \times H^{1}(\R)$ the Cauchy problem (decoupled system)

\begin{equation}\label{ecua153}
\left\{
\begin{aligned}
 \partial_{t} u
 + i \ (-\Delta)^{s} u
 - i \ \Delta u
  & = - i \ \alpha \ \tilde{v} \ \tilde{u}
  - i \ \ |\tilde{u}|^{2} \ \tilde{u},   ~~& x \in \R ~~ t > 0, 
 \\[5pt]
 \partial_{t} v
 -  \Delta v
  & = \beta \ (-\Delta)^{s/2} (|\tilde{u}|^{2})
  - (-\Delta)^{s/2} g(\tilde{v}),   ~~& x \in \R ~~ t > 0,
 \\[5pt]
 u(0,x) = u_{0}(x), & ~~ v(0,x) = v_{0}(x),
\end{aligned}
\right.
\end{equation}
admits a unique mild solution $(u,v) \in C \big( [0,T];H^{1}(\R) \big) \times C \big( [0,T];H^{1}(\R) \big)$.
\end{lemma}

\begin{proof}
First, we define for each $t \in (0,T)$
$$
   F(t):= - i \ \alpha \ \tilde{v}(t) \ \tilde{u}(t) - i \ \ |\tilde{u}|^{2}(t) \ \tilde{u}(t), 
   \quad
   G(t):= \beta(-\Delta)^{s/2} (|\tilde{u}|^{2})(t)
  - (-\Delta)^{s/2} g(\tilde{v})(t).
$$

\medskip
 \underline {Claim 1:} The complex value function $F \in C([0,T];L^{2}(\R))$.

\medskip
{Proof of Claim:} Indeed, for all $t \! \in [0,T]$, 
$|\tilde{u}|^{2}(t) \ \tilde{u}(t) \in H^{1}(\R)$,
$\tilde{u}(t) \ \tilde{v}(t) \in H^{1}(\R)$.
Then, for $h$ sufficiently small
$$
\begin{aligned}
  F(t+h) - F(t) &=
  i \ \alpha \ \Big(\tilde{v}(t) \ \tilde{u}(t) - \tilde{v}(t+h) \ \tilde{u}(t+h) \Big)
\\[5pt]
  &\quad + i \ \ \Big( |\tilde{u}|^{2}(t) \ \tilde{u}(t) - |\tilde{u}|^{2}(t+h) \ \tilde{u}(t+h) \Big)
=  i \ \alpha \ I_{1} + i \ I_{2},
\end{aligned}
$$
with obvious notation. A simple algebraic computation shows that 
$$
    \lim_{h \to 0} \|I_1\|_{L^{2}(\R)}= 0, \quad \text{and} 
    \quad 
    \lim_{h \to 0} \|I_2\|_{L^{2}(\R)}= 0,
$$
from which the claim is proved. 

\medskip
 \underline {Claim 2:} The real value function $G \in C([0,T];L^{2}(\R))$.

\medskip
{Proof of Claim:} We observe that $(-\Delta)^{s/2} (|\tilde{u}|^{2})(t) \in L^{2}(\R)$,
for each $t \in (0,T)$.
Also from the assumptions for the function $g$, that is 
$g \in C^{1}(\R)$, $g(0) = 0$ and $|g^{\prime}(v)| \leq M$, $(\forall v \in \R)$, it follows that 
 $(-\Delta)^{s/2} g(\tilde{v})(t) \in L^{2}(\R)$. 
Now, for $h$ sufficiently small, we have 
$$
\begin{aligned}
G(t+h) - G(t)&= \beta \ \Big( (-\Delta)^{s/2} (|\tilde{u}|^{2})(t+h)
- (-\Delta)^{s/2} (|\tilde{u}|^{2})(t) \Big)
\\[5pt]
&\quad - \Big( (-\Delta)^{s/2} g(\tilde{v})(t+h) - (-\Delta)^{s/2} g(\tilde{v})(t) \Big)
 = \beta \ J_{1} - J_{2},
\end{aligned}
$$
with obvious notation.  
Then, from \eqref{NORMHs} and the embedding theorem 
$$
\begin{aligned}
   \|J_{1}\|^{2}_{L^{2}(\R)}
& \leq \| \ |\tilde{u}|^{2}(t+h)
   - |\tilde{u}|^{2}(t)\|_{H^{s}(\R)}^{2} 
\leq \| \ |\tilde{u}|^{2}(t+h) - |\tilde{u}|^{2}(t)\|_{H^{1}(\R)}^{2}.
\end{aligned}
$$
Analogously, we have
$$
\begin{aligned}
   \|J_{2}\|^{2}_{L^{2}(\R)}
   & \leq \|g(\tilde{v})(t+h)
   - g(\tilde{v})(t)\|_{H^{s}(\R)}^{2}
    \leq \|g(\tilde{v})(t+h)
   - g(\tilde{v})(t)\|_{H^{1}(\R)}^{2}
    \\[5pt]
   & = \|g(\tilde{v})(t+h)
   - g(\tilde{v})(t)\|_{L^{2}(\R)}^{2} + \|\partial _x g(\tilde{v})(t+h)
   - \partial _x g(\tilde{v})(t)\|_{L^{2}(\R)}^{2}
    \\[5pt]
   & \leq M^{2} \Big( \int_{\R} |\tilde{v}(t+h,x)
   - \tilde{v}(t,x)|^{2} \ dx 
   +  \int_{\R} \big| \partial _x \tilde{v}(t+h,x) - \partial _x \tilde{v}(t,x) \big|^{2} \ dx \Big)
\\[5pt]   
    & \leq 2 M^{2} \ \|\tilde{v}(t+h) - \tilde{v}(t)\|_{H^{1}(\R)}^{2}.
\end{aligned}
$$
Then, passing to the limit as $h \to 0$, the claim is proved. 

\bigskip
Finally, since $F, G \in C([0,T];L^{2}(\R))$ applying Lemma 4.15 and Corollary 4.12 in \cite{Cazanave}, there exists a 
unique solution $(u,v) \in C \big( [0,T];H^{1}(\R) \big) \times C \big( [0,T];H^{1}(\R) \big)$ given by
\begin{equation}
\label{SOLDES}
\begin{aligned}
   u(t)&= U(t) \ u_{0} - i \int_{0}^{t}U(t-t^{\prime}) \ \big(\alpha \ \widetilde{v}(t^{\prime}) \ \tilde{u}(t^{\prime})
+ \ |\tilde{u}(t^{\prime})|^{2} \ \tilde{u}(t^{\prime})\big) \ dt^{\prime},
\\[5pt]
v(t)&= W(t) \ v_{0}
+ \beta \int_{0}^{t} W \big( t-t^{\prime} \big) \  \big( \beta \ (-\Delta)^{s/2} (|\tilde{u}(t^{\prime})|^{2})
- (-\Delta)^{s/2} g(\tilde{v})(t^{\prime}) \big) \ dt^{\prime},
\end{aligned}
\end{equation}
where $U(t)\equiv U_{\ve=1}(t)$,  $W(t)\equiv W_{\ve=1}(t)$ are given respectively by 
\eqref{GROUPS}, and \eqref{GROUPH}.
\end{proof}


\begin{proposition}
\label{sollocal}
Let $\frac{1}{2} < s < 1$, $g \in C^{1}(\R)$, $0 < m \leq g^{\prime}(\cdot) \leq M < \infty$, $(g(0)=0)$. 
Then, for any $(u^\ve_{0},v^\ve_{0}) \in H^{1}(\R) \times H^{1}(\R)$, there exists $T > 0$ such that, 
the Cauchy problem \eqref{ecuaI2} has a unique mild solution.
\end{proposition}

\begin{proof}
1. Hereupon, we denote by $X_T$ the Banach space $C \big( [0,T];H^{1}(\R) \big)$, 
where $T> 0$ is chosen a posteriori. 
For $R > 2 \max\{\|u^\ve_{0}\|_{H^{1}(\R)}, \|v^\ve_{0}\|_{H^{1}(\R)}\}$, we define
$$
    B_R^T:= \{ f \in X_T : \|f\|_{L^{\infty}(0,T;H^{1}(\R))} \leq R \}, 
$$
and the mapping $\Phi: B_R^T \times  B_R^T \to X_T \times X_T$,
$(\tilde{u},\tilde{v}) \mapsto (u^\ve,v^\ve) \equiv \Phi(\tilde{u},\tilde{v})$,
where $(u^\ve,v^\ve)$ is the unique mild solution of the Cauchy problem \eqref{ecua153} 
(for each $\ve> 0$ fixed). 
Then, from \eqref{SOLDES} we have for any $t \in [0,T]$
$$
\begin{aligned}
\Phi_{1}(\tilde{u},\tilde{v})\equiv
u^{\ve}(t)&= U_{\ve}(t) u^\ve_{0}
- i \int_{0}^{t} U_{\ve}(t-t^{\prime}) \ \big( \alpha \ \tilde{v}(t^{\prime}) \ \tilde{u}(t^{\prime})
+ \ |\tilde{u}(t^{\prime})|^{2} \ \tilde{u}(t^{\prime}) \big) \ dt^{\prime},
\\[5pt]
\Phi_{2}(\tilde{u},\tilde{v})\equiv
v^{\ve}(t)&= W_{\ve}(t) v^\ve_{0}
+ \int_{0}^{t} W_{\ve} \big( t-t^{\prime} \big) \ \big( \beta \ (-\Delta)^{s/2} (|\tilde{u}(t^{\prime})|^{2})
- (-\Delta)^{s/2} g(\tilde{v})(t^{\prime}) \big) \ dt^{\prime}.
\end{aligned}
$$

\bigskip
2. First, we show that 
$(\Phi_{1}(\tilde{u},\tilde{v}),\Phi_{2}(\tilde{u},\tilde{v})) \in B_{R}^{T} \times B_{R}^{T}$.
Indeed, since for each $t \in [0,T]$, $\|U_{\ve}(t) u^\ve_{0}\|_{H^{1}(\R)}= \|u^\ve_{0}\|_{H^{1}(\R)}$,
then
$$
    \|U_{\ve}(\cdot) u^\ve_{0}\|_{L^{\infty}(0,T;H^{1}(\R))} = \|u^\ve_{0}\|_{H^{1}(\R)}.
$$
Moreover, we have
$$
\begin{aligned}
\|\int_{0}^{t}& U_{\ve}(t-t^{\prime}) \ \big( \alpha \ \tilde{v}(t^{\prime}) \ \tilde{u}(t^{\prime})
+ \ |\tilde{u}(t^{\prime})|^{2} \ \tilde{u}(t^{\prime}) \big) \ dt^{\prime}\|_{H^{1}(\R)}
\\[5pt]
&\leq \int_{0}^{t} \|\big( \alpha \ \tilde{v}(t^{\prime}) \ \tilde{u}(t^{\prime})
+ \ |\tilde{u}(t^{\prime})|^{2} \ \tilde{u}(t^{\prime}) \big)\|_{H^{1}(\R)} \ dt^{\prime}
\\[5pt]
& \leq \int_{0}^{t} C \ \big( |\alpha| \ \|\tilde{v}(t^{\prime})\|_{H^{1}(\R)} \ \|\tilde{u}(t^{\prime})\|_{H^{1}(\R)}
+ \ \||\tilde{u}(t^{\prime})|^{2}\|_{H^{1}(\R)} \ \|\tilde{u}(t^{\prime})\|_{H^{1}(\R)} \big) \ dt^{\prime}
\\[5pt]
& = C \ |\alpha| \int_{0}^{t} \|\tilde{v}(t^{\prime})\|_{H^{1}(\R)} \ \|\tilde{u}(t^{\prime})\|_{H^{1}(\R)} \ dt^{\prime}
+ C \  \int_{0}^{t} \|\tilde{u}(t^{\prime})\|_{H^{1}(\R)}^{3} \ dt^{\prime}
\\[5pt]
& \leq 2 \max\{|\alpha|, R\} \ C \ R^{2} \ T,
\end{aligned}
$$
where we have used \eqref{algebra}. Consequently, for $T$ satisfying 
\begin{equation}
\label{ecua159}
 T< \frac{1}{4 \max\{|\alpha|, R\} \ C \ R}, 
\end{equation}
$$
\begin{aligned}
    \|\Phi_{1}(\tilde{u},\tilde{v})\|_{L^{\infty}(0,T;H^{1}(\R))}
    &\leq \|u^\ve_{0}\|_{H^{1}(\R)} + 2 \max\{|\alpha|, R\} \ C \ R^{2} \ T
    < \frac{R}{2} + \frac{R}{2}= R.
\end{aligned}
$$

Similarly, we estimate $\|\Phi_{2}(\tilde{u},\tilde{v})\|_{L^{\infty}(0,T;H^{1}(\R))}$. Applying \eqref{E03}, it follows that 
$$
\begin{aligned}
   \|\int_{0}^{t} & W_{\ve} \big( t-t^{\prime} \big) \ \big( \beta \ (-\Delta)^{s/2} (|\tilde{u}(t^{\prime})|^{2}) 
   - (-\Delta)^{s/2} g(\tilde{v})(t^{\prime}) \big) \ dt^{\prime}\|_{H^{1}(\R)}
\\[5pt]
& \leq \int_{0}^{t}\big(1 +  \frac{C}{(t-t^{\prime})^{1/2}}\big) \|\beta \ (-\Delta)^{s/2} (|\tilde{u}(t^{\prime})|^{2}) 
- (-\Delta)^{s/2} g(\tilde{v})(t^{\prime})\|_{L^{2}(\R)} \ dt^{\prime}
\\[5pt]
& \leq \int_{0}^{t} |\beta| \ \big(1 +  \frac{C}{(t-t^{\prime})^{1/2}}\big) \  \|(-\Delta)^{s/2} (|\tilde{u}(t^{\prime})|^{2})\|_{L^{2}(\R)} \ dt^{\prime}
\\[5pt]
& \leq \int_{0}^{t}\big(1 +  \frac{C}{(t-t^{\prime})^{1/2}}\big) \ \|(-\Delta)^{s/2} g(\tilde{v})(t^{\prime})\|_{L^{2}(\R)} \ dt^{\prime}= I_1 + I_2,
\end{aligned}
$$
with obvious notation. To follow, we have
$$
\begin{aligned}
I_1 &\leq \frac{1}{m_s} \int_{0}^{t} |\beta| \ \big(1 +  \frac{C}{(t-t^{\prime})^{1/2}}\big) \  \|\tilde{u}(t^{\prime})\|^{2}_{H^{s}(\R)} \ dt^{\prime}
\\[5pt]
&\leq \frac{|\beta|}{m_s} \int_{0}^{t} \ \big(1 +  \frac{C}{(t-t^{\prime})^{1/2}}\big)  \  \|\tilde{u}(t^{\prime})\|^{2}_{H^{1}(\R)} \ dt^{\prime}
< \frac{|\beta|}{m_s} \ R^{2} \ ( T + 2 C \sqrt{T}), 
\end{aligned}
$$
and 
$$
\begin{aligned}
I_2 &\leq  \frac{1}{m_s} \int_{0}^{t}\big(1 +  \frac{C}{(t-t^{\prime})^{1/2}}\big) \ \|g(\tilde{v})(t^{\prime})\|_{H^{s}(\R)} \ dt^{\prime}
\\[5pt]
&\leq  \frac{1}{m_s} \int_{0}^{t}\big(1 +  \frac{C}{(t-t^{\prime})^{1/2}}\big) \ 
\big(\|g(\tilde{v})(t^{\prime})\|_{L^{2}(\R)}^{2} + \|\partial_{x}g(\tilde{v})(t^{\prime})\|_{L^{2}(\R)}^{2} \big)^{1/2}  \ dt^{\prime}
\\[5pt]
&\leq \frac{M}{m_s} \int_{0}^{t}\big(1 +  \frac{C}{(t-t^{\prime})^{1/2}}\big) \ \|\tilde{v})(t^{\prime})\|_{H^{1}(\R)} \ dt^{\prime}
< \frac{M}{m_s}  \ R \ ( T + 2 C \sqrt{T}).
\end{aligned}
$$
Consequently, for $T$ satisfying 
\begin{equation}
\label{ecua1599}
 T< \min\{ \frac{m_s}{8 \max\{|\beta| R, M\}},  \frac{(m_s)^{2}}{64 C^{2} (\max\{|\beta| R, M\})^{2}} \}, 
\end{equation}
$$
\begin{aligned}
    \|\Phi_{2}(\tilde{u},\tilde{v})\|_{L^{\infty}(0,T;H^{1}(\R))}
    &\leq \|v^\ve_{0}\|_{H^{1}(\R)} + \frac{R}{m_s} (|\beta| R + M) (T + \sqrt{T})
    < \frac{R}{2} + \frac{R}{2}= R.
\end{aligned}
$$

\bigskip
3. Now, we show that $\Phi$ is a contraction on $B_{R}^{T} \times B_{R}^{T}$. Let 
$(\tilde{u}_i,\tilde{v}_i) \in B_R^T \times  B_R^T$, $(i= 1, 2)$, then we have 
\begin{equation}
\label{CONTR11}
\begin{aligned}
  \|\Phi_{1}(\tilde{u}_1,\tilde{v}_1) &- \Phi_{1}(\tilde{u}_2,\tilde{v}_2)\|_{H^{1}(\R)}
  \\[5pt]
 & \leq  |\alpha| \int_{0}^{t} \| U_{\ve}(t-t^{\prime}) \ \big( \tilde{v}_2(t^{\prime}) \ \tilde{u}_2(t^{\prime})
   - \tilde{v}_1(t^{\prime}) \ \tilde{u}_1(t^{\prime}) \big)\|_{H^{1}(\R)} \ dt^{\prime}
    \\[5pt]
   &  + \int_{0}^{t} \| U_{\ve}(t-t^{\prime}) \big( |\tilde{u}_2(t^{\prime})|^{2} \ \tilde{u}_2(t^{\prime})  
   - |\tilde{u}_1(t^{\prime})|^{2} \ \tilde{u}_1(t^{\prime})
   \big) \|_{H^{1}(\R)} \ dt^{\prime}
     \\[5pt]
 & \leq  |\alpha| \int_{0}^{t} \| \tilde{v}_2(t^{\prime}) \ \tilde{u}_2(t^{\prime})
   - \tilde{v}_1(t^{\prime}) \ \tilde{u}_1(t^{\prime}) \|_{H^{1}(\R)} \ dt^{\prime}
    \\[5pt]
   &  + \int_{0}^{t} \| \ |\tilde{u}_2(t^{\prime})|^{2} \ \tilde{u}_2(t^{\prime})  
   - |\tilde{u}_1(t^{\prime})|^{2} \ \tilde{u}_1(t^{\prime}) \|_{H^{1}(\R)} \ dt^{\prime} = |\alpha| \ J_{1} + J_{2}.
\end{aligned}
\end{equation}
Applying \eqref{algebra} we obtain 
\begin{equation}
\label{CONTR12}
\begin{aligned}
|\alpha| J_{1} &\leq \ C \ |\alpha|\int_{0}^{t} \|\tilde{v}_2(t^{\prime})\|_{H^{1}(\R)} \ \|\tilde{u}_2(t^{\prime})
- \tilde{u}_1(t^{\prime})\|_{H^{1}(\R)} \ dt^{\prime}
 \\[5pt]
& + \ C \ |\alpha|\int_{0}^{t} \|\tilde{u}_1(t^{\prime})\|_{H^{1}(\R)} \ \|\tilde{v}_2(t^{\prime})
- \tilde{v}_1(t^{\prime})\|_{H^{1}(\R)} \ dt^{\prime}
 \\[5pt]
& \leq C \ |\alpha| \ R \ T \ \big(\|\tilde{u}_2 - \tilde{u}_1\|_{L^{\infty}(0,T;H^{1}(\R))}
+ \|\tilde{v}_2 - \tilde{v}_1 \|_{L^{\infty}(0,T;H^{1}(\R))} \big).
\end{aligned}
\end{equation}
Similarly, we also have
\begin{equation}
\label{CONTR13}
\begin{aligned}
J_2 & \leq C \int_0^t \big( \| \tilde{u}_2(t^{\prime})\|^{2}_{H^{1}(\R)} \ \|\tilde{u}_2(t^{\prime})
- \tilde{u}_1(t^{\prime})\|_{H^{1}(\R)}
\\[5pt]
& \quad  \quad + \|\tilde{u}_1(t^{\prime})\|_{H^{1}(\R)} \ \| \ |\tilde{u}_2(t^{\prime})|^{2}
- |\tilde{u}_1(t^{\prime})|^{2}\|_{H^{1}(\R)} \big)  \ dt^{\prime}
 \\[5pt]
& \leq 3 C \ R^{2}  \int_0^t \|\tilde{u}_2(t^{\prime}) - \tilde{u}_1(t^{\prime})\|_{H^{1}(\R)} \ dt^{\prime}
\leq 3 C \ R^{2} \ T \ \|\tilde{u}_2
- \tilde{u}_1\|_{L^{\infty}(0,T;H^{1}(\R))}.
\end{aligned}
\end{equation}
Therefore, from \eqref{CONTR11}--\eqref{CONTR13}, it follows that 
$$
\begin{aligned}
  \|\Phi_{1}(\tilde{u}_1,&\tilde{v}_1) - \Phi_{1}(\tilde{u}_2,\tilde{v}_2)\|_{H^{1}(\R)}
  \\[5pt]
 & \leq C \ R \ \max \{|\alpha|, 3 R \} \ T  \ \big(\|\tilde{u}_1 - \tilde{u}_2\|_{L^{\infty}(0,T;H^{1}(\R))}
+ \|\tilde{v}_1 - \tilde{v}_2 \|_{L^{\infty}(0,T;H^{1}(\R))} \big).
 \end{aligned}
$$

To this end, we have
\begin{equation}\label{E011}
\begin{aligned}
  \|\Phi_{2}(\tilde{u}_1,\tilde{v}_1) & - \Phi_{2}(\tilde{u}_2,\tilde{v}_2)\|_{H^{1}(\R)}
\\[5pt]
& \leq \int_{0}^{t} |\beta| \ \|W_{\ve} \big( t-t^{\prime} \big) \big( (-\Delta)^{s/2} |\tilde{u}_1(t^{\prime})|^{2}
- (-\Delta)^{s/2} |\tilde{u}_2(t^{\prime})|^{2} \big) \|_{H^{1}(\R)} \ dt^{\prime}
\\[5pt]
& +  \int_{0}^{t} \|W_{\ve} \big( t-t^{\prime} \big)  \big( (-\Delta)^{s/2} g(\tilde{v}_2)(t^{\prime})
- (-\Delta)^{s/2} g(\tilde{v}_1)(t^{\prime}) \big) \|_{H^{1}(\R)} \ dt^{\prime}
\\[5pt]
& \leq \int_{0}^{t} |\beta| \ \big(1 +  \frac{C}{(t-t^{\prime})^{1/2}}\big) \ \| (-\Delta)^{s/2} \big( |\tilde{u}_1(t^{\prime})|^{2}
- |\tilde{u}_2(t^{\prime})|^{2} \big) \|_{L^{2}(\R)} \ dt^{\prime}
\\[5pt]
& +  \int_{0}^{t} \ \big(1 +  \frac{C}{(t-t^{\prime})^{1/2}}\big) \ \|(-\Delta)^{s/2} \big(g(\tilde{v}_2)(t^{\prime})
- g(\tilde{v}_1)(t^{\prime}) \big) \|_{L^{2}(\R)} \ dt^{\prime}
\\[5pt]
&= K_1 + K_2, 
\end{aligned}
\end{equation}
where we have used \eqref{E03}, and obvious notation. Applying \eqref{NORMHs}, we obtain
\begin{equation}
\label{E012}
\begin{aligned}
K_1 &\leq \frac{|\beta|}{m_s} \int_{0}^{t} \ \big(1 +  \frac{C}{(t-t^{\prime})^{1/2}}\big)  
\  \| \ |\tilde{u}_1(t^{\prime})|^{2} - \tilde{u}_2(t^{\prime})|^{2}\|_{H^{s}(\R)} \ dt^{\prime}
\\[5pt]
&\leq \frac{|\beta|}{m_s} \int_{0}^{t} \ \big(1 +  \frac{C}{(t-t^{\prime})^{1/2}}\big)  
\   \bigg( \|\tilde{u}_1(t^{\prime})\|_{H^{1}(\R)} \ \|\tilde{u}_1(t^{\prime}) 
- \tilde{u}_2(t^{\prime})\|_{H^{1}(\R)}
\\[5pt]
& \hspace{120pt} + \|\tilde{u}_2(t^{\prime})\|_{H^{1}(\R)} \ \|\tilde{u}_1(t^{\prime}) 
- \tilde{u}_2(t^{\prime})\|_{H^{1}(\R)} \bigg) \ dt^{\prime}
\\[5pt]
&\leq 2 R \  \frac{|\beta|}{m_s} \ (T + C \sqrt{T}) \ \|\tilde{u}_1 - \tilde{u}_2\|_{L^{\infty}(0,T;H^{1}(\R))},
\end{aligned}
\end{equation}
and 
\begin{equation}
\label{E0122}
\begin{aligned}
K_2 &\leq \frac{1}{m_s}  \int_{0}^{t} \ \big(1 +  \frac{C}{(t-t^{\prime})^{1/2}}\big)  
\  \| g(\tilde{v}_2)(t^{\prime})
- g(\tilde{v}_1)(t^{\prime}) \|_{H^{s}(\R)} \ dt^{\prime}
\\[5pt]
 &\leq \frac{1}{m_s}  \int_{0}^{t} \ \big(1 +  \frac{C}{(t-t^{\prime})^{1/2}}\big)  
   \bigg( \|g(\tilde{v}_2)(t^{\prime}) - g(\tilde{v}_1)(t^{\prime})\|_{L^{2}(\R)}^{2}
\\[5pt]
& \hspace{120pt} + \|\partial _x g(\tilde{v}_2)(t^{\prime}) 
- \partial _x g(\tilde{v}_1)(t^{\prime})\|_{L^{2}(\R)}^{2} \bigg)^{1/2} \ dt^{\prime}
\\[5pt]
&\leq   \frac{M}{m_s} \ (T + C \sqrt{T}) \ \|\tilde{v}_1 - \tilde{v}_2\|_{L^{\infty}(0,T;H^{1}(\R))}.
\end{aligned}
\end{equation}
Consequently, from \eqref{E011}--\eqref{E0122} we obtain
$$
\begin{aligned}
  \|\Phi_{2}(\tilde{u}_1,&\tilde{v}_1) - \Phi_{2}(\tilde{u}_2,\tilde{v}_2)\|_{H^{1}(\R)}
  \\[5pt]
 & \leq \frac{\max \{ 2 R \ |\beta|, M \}}{m_s} \ (T + \sqrt{T})  \ \big(\|\tilde{u}_1 - \tilde{u}_2\|_{L^{\infty}(0,T;H^{1}(\R))}
+ \|\tilde{v}_1 - \tilde{v}_2 \|_{L^{\infty}(0,T;H^{1}(\R))} \big).
 \end{aligned}
$$

\bigskip
4. Finally, from items (2) and (3) there exists a $T> 0$, sufficiently small, such that 
$\Phi: B_R^T \times  B_R^T \to B_R^T \times B_R^T$ is a (strict) contraction. Hence we 
can apply the Banach Fixed Point Theorem and obtain a unique (local in time) solution $(u^\ve, v^\ve)$ of
the Cauchy problem \eqref{ecuaI2}. 
\end{proof}

\subsection{A priori estimates}

For each $\ve> 0$, let $(u^\ve, v^\ve)$ be the unique solution for the Cauchy problem \eqref{ecuaI2}, 
and recall that, the sequences $\{u_0^\ve\}$ and $\{v_0^\ve\}$ are uniformly bounded in 
$H^1(\R)$ with respect to $\ve> 0$ fixed. 
\begin{lemma}[First estimate]
\label{First estimate} 
Let $\frac{1}{2} < s < 1$.
Then, for each $t \in (0,T)$
\begin{equation}
\label{ecuaI9}
    \frac{d}{dt} \int_{\R} |u^{\ve}(t,x)|^{2} \ dx= 0,
\end{equation}
\begin{equation}\label{ecuaI10}
\begin{aligned}
  \frac{d}{dt} \Big(& \int_{\R} |(-\Delta)^{s/2} u^{\ve}(t,x)|^{2} \ dx
  + \ve^{a} \int_{\R} |\partial_x u^{\ve}(t,x)|^{2} \ dx
  + \frac{1}{2} \int_{\R} |u^{\ve}(t,x)|^{4} \ dx
  \\[5pt]
  & + \alpha \int_{\R} v^{\ve}(t,x) \ |u^{\ve}(t,x)|^{2} \ dx \Big)
  = \alpha \ \beta \int_{\R} (-\Delta)^{s/2} (|u^{\ve}(t,x)|^{2}) \ |u^{\ve}(t,x)|^{2} \ dx
  \\[5pt]
  & - \alpha \int_{\R} |u^{\ve}(t,x)|^{2} \ (-\Delta)^{s/2} g_{\ve}(v^{\ve}(t,x)) \ dx
  - \alpha \ \ve^{b} \int_{\R} \partial_x |u^{\ve}(t,x)|^{2} \ \partial_x v^{\ve}(t,x) \ dx,
\end{aligned}
\end{equation}
\begin{equation}
\label{ecuaI11}
\begin{aligned}
& \frac{1}{2} \frac{d}{dt} \int_{\R} |v^{\ve}(t,x)|^{2} \ dx
    + \int_{\R} (-\Delta)^{s/2} g(v^{\ve})(t,x) \ v^{\ve}(t,x) \ dx 
  \\[5pt]
  & \quad + \ve^{b} \int_{\R} |\partial_x v^{\ve}(t,x)|^{2} \ dx = \beta \int_{\R} (-\Delta)^{s/2} (|u^{\ve}(t,x)|^{2}) \ v^{\ve}(t,x) \ dx.
\end{aligned}    
\end{equation}
\end{lemma}

\begin{proof}
1. First, by approximating the initial data in $H^1(\R)$ 
by functions in $C^\infty_c(\R)$, and a standard 
limit argument, we can assume that $(u^\ve, v^\ve)$
satisfies the Cauchy problem \eqref{ecuaI2}, (at least
almost everywhere), and we are allowed to 
make the computations below. Indeed, since 
$H^s(\R)$ is an algebra for any $s> 1/2$,
we may follow the same 
strategy developed in the previous section, and 
for $0< T'< T$, we obtain  
$(u^\ve, v^\ve) \in \big(C([0,T']; H^k(\R)) \cap C^1([0,T']; H^{k-2}(\R))\big)^2$,
for each integer $k> 2$. 

\medskip
2. To follow, multiplying equation \eqref{ecuaI2}$_{1}$ by $\overline{u^{\ve}}(t,x)$
and integrating in $\R$, we have
$$
\begin{aligned}
& i \int_{\R} \partial_{t} u^{\ve}(t,x) \ \overline{u^{\ve}}(t,x) \ dx - \int_{\R} |(-\Delta)^{s/2} u^{\ve}(t,x)|^{2} \ dx
- \ve^{a} \int_{\R} |\partial_x u^{\ve}(t,x)|^{2} \ dx
  \\[5pt]
& = \alpha \int_{\R} v^{\ve}(t,x) \ |u^{\ve}(t,x)|^{2} \ dx
+ \int_{\R} |u^{\ve}(t,x)|^{4} \ dx.
\end{aligned}  
$$
Therefore, taking the imaginary part of the above equation, we obtain
$$
\frac{1}{2} \frac{d}{dt} \int_{\R} |u^{\ve}(t,x)|^{2} \ dx
= \textrm{Re} \int_{\R} \partial_{t}u^{\ve}(t,x) \ \overline{u^{\ve}}(t,x) \ dx = 0.
$$

\medskip
3. Now, let us multiply equation \eqref{ecuaI2}$_{1}$ by $\partial_{t} \overline{u^{\ve}}(t,x)$, 
and integrate in $\R$ to obtain
$$
\begin{aligned}
& i \int_{\R} \partial_{t} u^{\ve}(t,x) \ \partial_{t} \overline{u^{\ve}}(t,x) \ dx
 - \int_{\R} (-\Delta)^{s} u^{\ve}(t,x) \ \partial_{t}\overline{u^{\ve}}(t,x) \ dx
\\[5pt] 
& \quad +  \ve^{a} \int_{\R} \Delta u^{\ve}(t,x) \ \partial_{t}\overline{u^{\ve}}(t,x) \ dx
  \\[5pt]
& \quad = \alpha \int_{\R} v^{\ve}(t,x) \ u^{\ve}(t,x) \ \partial_{t}\overline{u^{\ve}}(t,x) \ dx + \int_{\R} |u^{\ve}(t,x)|^{2} \ u^{\ve}(t,x) \ \partial_{t}\overline{u^{\ve}}(t,x) \ dx.
\end{aligned}     
$$
Then, writing $u^\ve= u_1^\ve + i u_2^\ve$ and integrating by parts, it follows that 
$$
\begin{aligned}
& i \int_{\R} |\partial_{t} u^{\ve}(t,x)|^{2} \ dx
- \int_{\R} (-\Delta)^{s/2}u^{\ve}(t,x) \ \partial_{t} \overline{(-\Delta)^{s/2} u^{\ve}}(t,x) \ dx
\\[5pt]
& \quad \quad - \ve^{a} \int_{\R} \partial_x u^{\ve}(t,x) \ \partial_{t}\overline{\partial_x u^{\ve}}(t,x) \ dx
  \\[5pt]
& = \alpha \int_{\R} v^{\ve}(t,x) \ (u_{1}^{\ve}(t,x)\partial_{t} u_{1}^{\ve}(t,x)+u_{2}^{\ve}(t,x)\partial_{t} u_{2}^{\ve}(t,x))  dx 
   \\[5pt]
& \quad \quad  + i \ \alpha \int_{\R} v^{\ve}(t,x) \ ( u_{2}^{\ve}(t,x)\partial_{t} u_{1}^{\ve}(t,x)-u_{1}^{\ve}(t,x)\partial_{t} u_{2}^{\ve}(t,x)) \ dx 
\\[5pt]
& \quad \quad + \frac{1}{2} \int_{\R} (u^{\ve}(t,x))^{2} \ \partial_{t}(\overline{(u^{\ve})^{2}}(t,x)) \ dx.
\end{aligned}
$$
Taking the real part we have
\begin{equation}\label{ecuaI12}
\begin{aligned}
&  \frac{d}{dt} \bigg[ \int_{\R} |(-\Delta)^{s/2} u^{\ve}(t,x)|^{2} \ dx + \ve^{a} \int_{\R} |\partial_x u^{\ve}(t,x)|^{2} \ dx
  + \frac{1}{2} \int_{\R} |u^{\ve}(t,x)|^{4} \ dx
  \\[5pt]
&  + \alpha \int_{\R} v^{\ve}(t,x) \ |u^{\ve}(t,x)|^{2} \ dx \bigg] = \alpha \int_{\R} |u^{\ve}(t,x)|^{2} \ \partial_{t} v^{\ve}(t,x) \ dx.
\end{aligned}
\end{equation}
The right-hand side of the above equation is computed by multiplying \eqref{ecuaI2}$_{2}$ by $\alpha |u^{\ve}(t,x)|^{2}$ and integrating in $\R$, 
that is to say
$$
\begin{aligned}
& \alpha \int_{\R} |u^{\ve}(t,x)|^{2} \ \partial_{t} v^{\ve}(t,x) \ dx = \alpha \ \beta \int_{\R} (-\Delta)^{s/2} (|u^{\ve}|^{2})(t,x) \ |u^{\ve}(t,x)|^{2} \ dx
 \\[5pt]
& - \alpha \int_{\R} |u^{\ve}(t,x)|^{2} \ (-\Delta)^{s/2} g_{\ve}(v^{\ve})(t,x) \ dx - \alpha \ \ve^{b} \int_{\R} \partial_x |u^{\ve}(t,x)|^{2} \ \partial_x v^{\ve}(t,x) \ dx,
\end{aligned}
$$
and replacing it in \eqref{ecuaI12}, we obtain
$$
\begin{aligned}
&  \frac{d}{dt} \bigg[ \int_{\R} |(-\Delta)^{s/2} u^{\ve}(t,x)|^{2} \ dx + \ve^{a} \int_{\R} |\partial_x u^{\ve}(t,x)|^{2} \ dx
  + \frac{1}{2} \int_{\R}|u^{\ve}(t,x)|^{4} \ dx
  \\[5pt]
&  + \alpha \int_{\R}v^{\ve}(t,x) \ |u^{\ve}(t,x)|^{2} \ dx \bigg] = \alpha \ \beta \int_{\R} (-\Delta)^{s/2} (|u^{\ve}|^{2})(t,x) \ |u^{\ve}(t,x)|^{2} \ dx
 \\[5pt]
&  - \alpha \int_{\R}|u^{\ve}(t,x)|^{2} \ (-\Delta)^{s/2} g_{\ve}(v^{\ve})(t,x) \ dx
  - \alpha \ \ve^{b} \int_{\R}\partial_x |u^{\ve}|^{2}(t,x) \ \partial_x v^{\ve}(t,x) \ dx.
\end{aligned}
$$

\medskip
3. Finally, equation \eqref{ecuaI11} follows directly by multiplying \eqref{ecuaI2}$_{2}$ 
by $v^{\ve}(t,x)$ and integrating in $\R$. Indeed, we have
$$
\begin{aligned}
& \frac{1}{2} \frac{d}{dt} \int_{\R} |v^{\ve}(t,x)|^{2} \ dx
+ \int_{\R} (-\Delta)^{s/2} g_{\ve}(v^{\ve})(t,x) \ v^{\ve}(t,x) \ dx
+ \ve^{b} \int_{\R} |\partial_x v^{\ve}(t,x)|^{2} \ dx
 \\[5pt]
& \quad = \beta \int_{\R} (-\Delta)^{s/2} (|u^{\ve}|^{2})(t,x) \ v^{\ve}(t,x) \ dx.
\end{aligned}
$$
\end{proof}

Now we pass to the second estimate.  
\begin{theorem}[Second estimate]
\label{Second estimate}
Let $\frac{1}{2} < s < 1$, and $g \in C^{1}(\R)$ satisfying
$$
    0 < \ve \leq g_{\ve}^{\prime}(\cdot)  \leq M.
$$    
Then for any $T>0$, there exist $\alpha_0> 0$ and $E_0> 0$, such that, for each $t \in (0,T)$
\begin{equation}
\label{ecuaI13}
\int_{\R} |(-\Delta)^{s/2} u^{\ve}(t,x)|^{2} \ dx
+ \ve^{4} \! \! \int_{\R} |\partial_x u^{\ve}(t,x)|^{2} \ dx + \frac{1}{4} \int_{\R} |u^{\ve}(t,x)|^{4} \ dx
\leq h(t),
\end{equation}
\begin{equation}
\label{ecuaI15FINAL}
\int_{\R} |v^{\ve}(t,x)|^{2} \ dx 
\leq e^{T}  \|v^{\ve}_{0}\|^{2}_{L^{2}(\R)}  
+ \frac{16 |\beta|^{2} e^{T}}{\pi(2s-1)} \|u^{\ve}_{0}\|_{L^{2}(\R)}^{2-\frac{1}{s}}\int_{0}^{t} h(\tau)^{1+\frac{1}{2s}}  \ d\tau \equiv H(t), \hspace{10pt}
\end{equation}
\begin{equation}
\label{ecuaI15}
\begin{aligned}
  &C_{1,s}^{-1} \int_{0}^{t} \| \ve^{1/2} \, (-\Delta)^{s/4} v^{\ve}(\tau)\|_{L^{2}(\R)}^{2} \ d\tau
 +  \int_{0}^{t} \| \ve^{7/2} \, \nabla v^{\ve}(\tau)\|^{2}_{L^{2}(\R)} \ d\tau \hspace{60pt}
    \\[5pt]
  &\leq \frac{\|v^{\ve}_{0}\|^{2}_{L^{2}(\R)}}{2} 
  + \frac{8 \,  |\beta|^{2}}{\pi(2s-1)} \|u^{\ve}_{0}\|_{L^{2}(\R)}^{2-\frac{1}{s}}  \int_{0}^{t} h(\tau)^{2+\frac{1}{s}} \ d\tau
  + \frac{1}{2} \int_0^t H^2(\tau) \, d\tau, 
\end{aligned}
\end{equation} 
for $|\alpha| \leq \alpha_0$ or 
$\|u_{0}\|_{L^{2}(\R)} \leq E_0$, where $h$ is a continuous positive function (independent of $\ve$).
\end{theorem}

\begin{proof}
1. First, from Proposition \ref{BWV11}
$$
    \int_{\R} (-\Delta)^{s/2} g_{\ve}(v^{\ve})(t,x) \ v^{\ve}(t,x) \ dx \geq \ve \ C_{1,s}^{-1} \ \|(-\Delta)^{s/4}v^{\ve}(t)\|_{L^{2}(\R)}^{2}.
$$
From the above inequality and equation \eqref{ecuaI11}, it follows that
$$
\begin{aligned}
&  \frac{1}{2} \frac{d}{dt} \int_{\R} |v^{\ve}(t,x)|^{2} \ dx
  + \ve \ C_{1,s}^{-1} \int_{\R} |(-\Delta)^{s/4} v^{\ve}(t,x)|^{2} \ dx
  + \ve^{b} \int_{\R} |\partial_x v^{\ve}(t,x)|^{2} \ dx
 \\[5pt]
&  \leq \beta \int_{\R} (-\Delta)^{s/4} |u^{\ve}(t,x)|^{2} \ (-\Delta)^{s/4}v^{\ve}(t,x) \ dx
 \\[5pt]
& \leq \frac{C_{1,s} \ \beta^{2}}{2\ve} \int_{\R} |(-\Delta)^{s/4} |u^{\ve}|^{2}(t,x)|^{2} \ dx
+ \frac{\ve \ C_{1,s}^{-1}}{2} \int_{\R} |(-\Delta)^{s/4} v^{\ve}(t,x)|^{2} \ dx,
\end{aligned}
$$
where we have used Young's inequality. 
Then, integrating from $0$ to $t> 0$,
\begin{equation}\label{ecuaI14}
\begin{aligned}
    & \int_{\R} |v^{\ve}(t,x)|^{2} \ dx
   + \ve \ C_{1,s}^{-1} \int_{0}^{t} \|(-\Delta)^{s/4} v^{\ve}(\tau)\|_{L^{2}(\R)}^{2} \ d\tau + 2 \ \ve^{b} \int_{0}^{t} \|\partial_x v^{\ve}(\tau)\|^{2}_{L^{2}(\R)} \ d\tau 
 \\[5pt]
   & \leq \int_{\R} |v^{\ve}_{0}(x)|^{2} \ dx
  + \frac{C_{1,s} \ \beta^{2}}{\ve} \int_{0}^{t}\|(-\Delta)^{s/4} |u^{\ve}|^{2}(\tau)\|_{L^{2}(\R)}^{2} \ d\tau.
\end{aligned}
\end{equation}

\medskip
2. Now, applying Proposition \ref{product} and 
equation \eqref{ecuaI9}, we have
\begin{equation}
\label{new}
\begin{aligned}
\|(-\Delta)^{s/2} |u^{\ve}|^{2}(t)\|_{L^{2}(\R)}
&\leq 2 \|u^{\ve}(t)\|_{L^{\infty}(\R)} \ \|(-\Delta)^{s/2} u^{\ve}(t)\|_{L^{2}(\R)}
\\[5pt]
&\leq \frac{4}{\sqrt{\pi(2s-1)}} \|u^{\ve}_{0}\|_{L^{2}(\R)}^{1-\frac{1}{2s}} \;  \|(-\Delta)^{s/2} u^{\ve}(t)\|_{L^{2}(\R)}^{1+\frac{1}{2s}}.
\end{aligned}
\end{equation}
Then, we obtain from \eqref{ecuaI14} and \eqref{new}
\begin{equation}
\label{ecuaI155}
\begin{aligned}
  \int_{\R} |v^{\ve}(t,x)|^{2} \ dx
  +  \ve \ C_{1,s}^{-1} \int_{0}^{t} \|(-\Delta)^{s/4} v^{\ve}(\tau)\|_{L^{2}(\R)}^{2} \ d\tau
  + 2 \ \ve^{b} \int_{0}^{t} \|\partial_x v^{\ve}(\tau)\|^{2}_{L^{2}(\R)} \ d\tau
    \\[5pt]
  \leq \|v^{\ve}_{0}\|^{2}_{L^{2}(\R)}
  + \frac{16 C_{1,s} \ \beta^{2}}{\ve \pi(2s-1)} \|u^{\ve}_{0}\|_{L^{2}(\R)}^{2-\frac{1}{s}}\int_{0}^{t} \|(-\Delta)^{s/2} u^{\ve}(\tau)\|_{L^{2}(\R)}^{2+\frac{1}{s}} \ d\tau.
\end{aligned}
\end{equation}
Similarly, we obtain 
$$
\begin{aligned}
& \int_{\R} |v^{\ve}(t,x)|^{2} \ dx 
  + 2 \ve \ C_{1,s}^{-1} \!\! \int_0^t \!\!\! \int_{\R} |(-\Delta)^{s/4} v^{\ve}(t,x)|^{2} \ dx d\tau
  + 2 \ve^{b} \!\! \int_0^t \!\!\!   \int_{\R} |\partial_x v^{\ve}(t,x)|^{2} \ dx d\tau
\\[5pt]
&\leq \|v^{\ve}_{0}\|^{2}_{L^{2}(\R)} + |\beta|^{2} \int_0^t \!\!\!  \int_{\R} |(-\Delta)^{s/2}| u^{\ve}(t,x)|^{2}|^{2} dx
+ \int_0^t \!\!\! \int_{\R} |v^{\ve}(t,x)|^{2} dx,
\end{aligned}
$$
and applying Gronwall's Lemma
$$
\int_{\R} |v^{\ve}(t,x)|^{2} \ dx 
\leq e^{T}  \|v^{\ve}_{0}\|^{2}_{L^{2}(\R)} + |\beta|^{2} e^{T} \int_{0}^{t} \|(-\Delta)^{s/2} |u^{\ve}(\tau)|^{2}\|^{2}_{L^{2}(\R)} d\tau.
$$
Therefore, from the above inequality and \eqref{new}, we have 
\begin{equation}
\label{ecuaD014}
\int_{\R} |v^{\ve}(t,x)|^{2} \ dx 
\leq e^{T}  \|v^{\ve}_{0}\|^{2}_{L^{2}(\R)}  
+ \frac{16 |\beta|^{2} e^{T}}{\pi(2s-1)} \|u^{\ve}_{0}\|_{L^{2}(\R)}^{2-\frac{1}{s}}\int_{0}^{t} \|(-\Delta)^{s/2} u^{\ve}(\tau)\|_{L^{2}(\R)}^{2+\frac{1}{s}} \ d\tau, 
\end{equation}
and
\begin{equation}
\label{new2}
\begin{aligned}
& \ve \ C_{1,s}^{-1} \!\! \int_0^t \!\!\! \int_{\R} |(-\Delta)^{s/4} v^{\ve}(t,x)|^{2} \ dx d\tau
  + \ve^{b} \!\! \int_0^t \!\!\!   \int_{\R} |\partial_x v^{\ve}(t,x)|^{2} \ dx d\tau \leq \frac{\|v^{\ve}_{0}\|^{2}_{L^{2}(\R)}}{2} 
\\[5pt]
  &+ \frac{8 \,  |\beta|^{2}}{\pi(2s-1)} \|u^{\ve}_{0}\|_{L^{2}(\R)}^{2-\frac{1}{s}}  \int_{0}^{t} \|(-\Delta)^{s/2} u^{\ve}(\tau)\|_{L^{2}(\R)}^{2+\frac{1}{s}}  \ d\tau
  + \frac{1}{2} \int_0^t \|v^\ve(\tau)\|^2_{L^2(\R)} \, d\tau. 
\end{aligned}
\end{equation}

\bigskip
Now, from equations \eqref{ecuaI10} it follows that 
\begin{equation}
\label{ecuaI16}
\begin{aligned}
   &  \frac{d}{dt} \bigg[ \int_{\R} |(-\Delta)^{s/2 }u^{\ve}(t,x)|^{2} \ dx
+ \ve^{a} \int_{\R} |\partial_x u^{\ve}(t,x)|^{2} \ dx
\\[5pt]
 &\quad + \frac{1}{2} \int_{\R} |u^{\ve}(t,x)|^{4} \ dx
+ \alpha \int_{\R} v^{\ve}(t,x) \ |u^{\ve}(t,x)|^{2} \ dx \bigg]
 \\[5pt]
& \quad \leq  |\alpha| \int_{\R} | (-\Delta)^{s/2} |u^{\ve}|^{2}(t,x) \  g_{\ve}(v^{\ve})(t,x) | dx
\\[5pt]
& \quad +  |\alpha| \ |\beta| \int_{\R} | (-\Delta)^{s/2} (|u^{\ve}|^{2})(t,x) | \ |u^{\ve}(t,x)|^{2} dx
\\[5pt]
& \quad + |\alpha| \ \ve^{b} \int_{\R} |\partial_x |u^{\ve}(t,x)|^{2} \ \partial_x v^{\ve}(t,x) | dx =: |\alpha| E + |\alpha| \ |\beta| F +  |\alpha| \ \ve^{b} G,
\end{aligned}
\end{equation}
with obvious notation. Again, from 
Proposition \ref{product} and equation \eqref{ecuaI9}, we may write:
$$
\begin{aligned} (i) \; 
   E 
  & \leq \| g_{\ve}(v^{\ve})(t)\|_{L^{2}(\R)} \ \|(-\Delta)^{s/2} |u^{\ve}|^{2}(t)\|_{L^{2}(\R)}
  \\[5pt]
   & \leq \|g^{\prime}_{\ve}\|_{L^{\infty}(\R)} \ \| v^{\ve}(t)\|_{L^{2}(\R)} \ 2 \ \|u^{\ve}(t)\|_{L^{\infty}(\R)} \ \|(-\Delta)^{s/2} u^{\ve}(t)\|_{L^{2}(\R)}
    \\[5pt]
   & \leq \frac{4}{\sqrt{\pi(2s-1)}} \ \|g^{\prime}_{\ve}\|_{L^{\infty}(\R)} \ \|u^{\ve}_{0}\|_{L^{2}(\R)}^{1-\frac{1}{2s}} 
   \ \| v^{\ve}(t)\|_{L^{2}(\R)} \ \|(-\Delta)^{s/2} u^{\ve}(t)\|_{L^{2}(\R)}^{1+\frac{1}{2s}}.
\end{aligned}
$$
$$
\begin{aligned}  (ii) \; 
  F 
  & \leq \|u^{\ve}(t)\|_{L^{\infty}(\R)} \int_{\R} |u^{\ve}(t,x)| \ |(-\Delta)^{s/2} (|u^{\ve}|^{2})(t,x)| \ dx
   \\[5pt]
   & \leq \|u^{\ve}(t)\|_{L^{\infty}(\R)} \ \|u^{\ve}(t)\|_{L^{2}(\R))} \ \|(-\Delta)^{s/2} |u^{\ve}|^{2}(t)\|_{L^{2}(\R))}
    \\[5pt]
   & \leq \frac{4}{\sqrt{\pi(2s-1)}} \ \|u^{\ve}_{0}\|_{L^{2}(\R)}^{2-\frac{1}{2s}} \ \|(-\Delta)^{s/2} u^{\ve}(t)\|_{L^{2}(\R)}^{\frac{1}{2s}} \ \|u^{\ve}(t)\|_{L^{\infty}(\R)} \ \|(-\Delta)^{s/2} u^{\ve}(t)\|_{L^{2}(\R)}
    \\[5pt]
   & 
   \leq \frac{8}{\pi(2s-1)} \ \|u^{\ve}_{0}\|_{L^{2}(\R)}^{3-\frac{1}{s}} \ \|(-\Delta)^{s/2} u^{\ve}(t)\|_{L^{2}(\R)}^{1+\frac{1}{s}}.
\end{aligned}
$$
\medskip
$$
\begin{aligned} (iii) \; 
  G 
   & \leq \int_{\R} \bigg( |\partial_x u^{\ve}(t,x)| \ |\overline{u^{\ve}}(t,x)| + |u^{\ve}(t,x)| \ |\partial_x \overline{u^{\ve}}(t,x)| \bigg) \ |\partial_x v^{\ve}(t,x)| \ dx
   \hspace{70pt}
   \\[5pt]
  & \leq 2 \|u^{\ve}(t)\|_{L^{\infty}(\R)} \int_{\R} |\partial_x u^{\ve}(t,x)| \ |\partial_x v^{\ve}(t,x)| \ dx
   \\[5pt]
   & \leq 2 \|u^{\ve}(t)\|_{L^{\infty}(\R)} \ \|\partial_x v^{\ve}(t)\|_{L^{2}(\R)} \ \|\partial_x u^{\ve}(t)\|_{L^{2}(\R)}
    \\[5pt]
   & \leq \frac{4}{\sqrt{\pi}} \ \|u^{\ve}(t)\|_{L^{2}(\R)}^{1/2} \ \|\partial_x u^{\ve}(t)\|_{L^{2}(\R)}^{1/2} \ \|\partial_x v^{\ve}(t)\|_{L^{2}(\R)} \ \|\partial_x u^{\ve}(t)\|_{L^{2}(\R)}
    \\[5pt]
   & \leq \frac{4}{\sqrt{\pi}} \ \|u^{\ve}_{0}\|_{L^{2}(\R)}^{1/2} \ \|\partial_x v^{\ve}(t)\|_{L^{2}(\R)} \ \|\partial_x u^{\ve}(t)\|_{L^{2}(\R)}^{3/2}.
\end{aligned}
$$
Replacing in equation \eqref{ecuaI16}
$$
\begin{aligned}
&  \frac{d}{dt} \Big( \int_{\R} |(-\Delta)^{s/2} u^{\ve}(t,x)|^{2} \ dx
+ \ve^{a} \int_{\R} |\partial_x u^{\ve}(t,x)|^{2} \ dx
\\[5pt]
& \quad + \frac{1}{2} \int_{\R} |u^{\ve}(t,x)|^{4} \ dx
+ \alpha \int_{\R} v^{\ve}(t,x) \ |u^{\ve}(t,x)|^{2} \ dx \Big)
 \\[5pt]
& \quad \leq |\alpha| \frac{4}{\sqrt{\pi(2s-1)}} \ \|g^{\prime}_{\ve}\|_{L^{\infty}(\R)} \ \|u^{\ve}_{0}\|_{L^{2}(\R)}^{1-\frac{1}{2s}} \ \| v^{\ve}(t)\|_{L^{2}(\R)}
   \|(-\Delta)^{s/2} u^{\ve}(t)\|_{L^{2}(\R)}^{1+\frac{1}{2s}}
\\[5pt]
& \quad + |\alpha| \ |\beta| \frac{8}{\pi(2s-1)} \ \|u^{\ve}_{0}\|_{L^{2}(\R)}^{3-\frac{1}{s}} \ \|(-\Delta)^{s/2} u^{\ve}(t)\|_{L^{2}(\R)}^{1+\frac{1}{s}}
\\[5pt]
& \quad + \frac{4}{\sqrt{\pi}} |\alpha| \ \ve^{b} \ \|u^{\ve}_{0}\|_{L^{2}(\R)}^{1/2} \ \|\partial_x v^{\ve}(t)\|_{L^{2}(\R)} \ \|\partial_x u^{\ve}(t)\|_{L^{2}(\R)}^{3/2}, 
\end{aligned}
$$
and integrating from $0$ to $t> 0$
$$
\begin{aligned}
&  \int_{\R} |(-\Delta)^{s/2} u^{\ve}(t,x)|^{2} \ dx
+ \ve^{a} \int_{\R} |\partial_x u^{\ve}(t,x)|^{2} \ dx
+ \frac{1}{2} \int_{\R} |u^{\ve}(t,x)|^{4} \ dx 
 \\[5pt]
& \quad \leq \|(-\Delta)^{s/2} u^{\ve}_{0}\|_{L^{2}(\R)}^{2}
+ \ve^{a} \ \|\partial_x u^{\ve}_{0}\|_{L^{2}(\R)}^{2} +\frac{1}{2} \|u^{\ve}_{0}\|_{L^{4}(\R)}^{4} + \int_{\R} |v^{\ve}_{0}(x)| \ |u^{\ve}_{0}(x)|^{2} \ dx
\\[5pt]
& \quad + \frac{4 \ |\alpha|}{\sqrt{\pi(2s-1)}} \ \|g^{\prime}_{\ve}\|_{L^{\infty}(\R)} \ \|u^{\ve}_{0}\|_{L^{2}(\R)}^{1-\frac{1}{2s}} \int_{0}^{t} \| v^{\ve}(\tau)\|_{L^{2}(\R)} \ \|(-\Delta)^{s/2} u^{\ve}(\tau)\|_{L^{2}(\R)}^{1+\frac{1}{2s}} \ d\tau
\\[5pt]
& \quad + \frac{8 \ |\alpha| \ |\beta|}{\pi(2s-1)} \ \|u^{\ve}_{0}\|_{L^{2}(\R)}^{3-\frac{1}{s}} \int_{0}^{t}\|(-\Delta)^{s/2} u^{\ve}(\tau)\|_{L^{2}(\R)}^{1+\frac{1}{s}} \ d\tau
\\[5pt]
& \quad + \frac{4}{\sqrt{\pi}} |\alpha| \ \ve^{b} \ \|u^{\ve}_{0}\|_{L^{2}(\R)}^{1/2} \int_{0}^{t} \|\partial_x v^{\ve}(\tau)\|_{L^{2}(\R)} \ \|\partial_x u^{\ve}(\tau)\|_{L^{2}(\R)}^{3/2} \ d\tau
\\[5pt]
& \quad + |\alpha| \int_{\R} |v^{\ve}(t,x)| \ |u^{\ve}(t,x)|^{2} \ dx.
\end{aligned}
$$
Then, we have
$$
\begin{aligned}
&  \int_{\R} |(-\Delta)^{s/2} u^{\ve}(t,x)|^{2} \ dx
+ \ve^{a} \int_{\R} |\partial_x u^{\ve}(t,x)|^{2} \ dx
+ \frac{1}{2} \int_{\R} |u^{\ve}(t,x)|^{4} \ dx 
 \\[5pt]
& \leq \|(-\Delta)^{s/2}u^{\ve}_{0}\|_{L^{2}(\R)}^{2}
+ \ve^{a} \ \|\partial_x u^{\ve}_{0}\|_{L^{2}(\R)}^{2} +\frac{1}{2} \|u^{\ve}_{0}\|_{L^{4}(\R)}^{4} + \|u^{\ve}_{0}\|_{L^{\infty}(\R)}\int_{\R} |v^{\ve}_{0}(x)| \ |u^{\ve}_{0}(x)| \ dx
\\[5pt]
& + \frac{4 \ |\alpha|}{\sqrt{\pi(2s-1)}} \ \|g^{\prime}_{\ve}\|_{L^{\infty}(\R)} \ \|u^{\ve}_{0}\|_{L^{2}(\R)}^{1-\frac{1}{2s}} 
\int_{0}^{t} \| v^{\ve}(\tau)\|_{L^{2}(\R)} \ \|(-\Delta)^{s/2} u^{\ve}(\tau)\|_{L^{2}(\R)}^{1+\frac{1}{2s}} \ d\tau
  \\[5pt]
& + \frac{8 \ |\alpha| \ |\beta|}{\pi(2s-1)} \ \|u^{\ve}_{0}\|_{L^{2}(\R)}^{3-\frac{1}{s}} \int_{0}^{t} \|(-\Delta)^{s/2} u^{\ve}(\tau)\|_{L^{2}(\R)}^{1+\frac{1}{s}} \ d\tau
\\[5pt]
& + \frac{4}{\sqrt{\pi}} |\alpha| \ \ve^{b} \ \|u^{\ve}_{0}\|_{L^{2}(\R)}^{1/2} \int_{0}^{t} \|\partial_x v^{\ve}(\tau)\|_{L^{2}(\R)} \ \|\partial_x u^{\ve}(\tau)\|_{L^{2}(\R)}^{3/2} \ d\tau
\\[5pt]
& + \int_{\R} \big(\sqrt{2} \ |\alpha| \ |v^{\ve}(t,x)| \big) \ \big( \frac{|u^{\ve}(t,x)|^{2}}{\sqrt{2}} \big) \ dx,
\end{aligned}
$$
from which follows that
\begin{equation}
\label{auxiliar}
\begin{aligned}
&  \int_{\R} |(-\Delta)^{s/2} u^{\ve}(t,x)|^{2} \ dx
+ \ve^{a} \int_{\R} |\partial_x u^{\ve}(t,x)|^{2} \ dx +\frac{1}{4} \int_{\R} |u^{\ve}(t,x)|^{4} \ dx 
  \\[5pt]
&  \quad \leq \|(-\Delta)^{s/2}u^{\ve}_{0}\|_{L^{2}(\R)}^{2}
+ \ve^{a} \ \|\partial_x u^{\ve}_{0}\|_{L^{2}(\R)}^{2}+\frac{1}{2} \|u^{\ve}_{0}\|_{L^{4}(\R)}^{4} + \|u^{\ve}_{0}\|_{L^{\infty}(\R)} \|v^{\ve}_{0} \|_{L^{2}(\R)} \|u^{\ve}_{0}\|_{L^{2}(\R)}
  \\[5pt]
& \quad  + \frac{4 \ |\alpha|}{\sqrt{\pi(2s-1)}} \ \|g^{\prime}_{\ve}\|_{L^{\infty}(\R)} \ \|u^{\ve}_{0}\|_{L^{2}(\R)}^{1-\frac{1}{2s}}\int_{0}^{t} \| v^{\ve}(\tau)\|_{L^{2}(\R)} \ \|(-\Delta)^{s/2} u^{\ve}(\tau)\|_{L^{2}(\R)}^{1+\frac{1}{2s}} \ d\tau
  \\[5pt]
& \quad + \frac{8 \ |\alpha| \ |\beta|}{\pi(2s-1)} \ \|u^{\ve}_{0}\|_{L^{2}(\R)}^{3-\frac{1}{s}} \int_{0}^{t} \|(-\Delta)^{s/2}u^{\ve}(\tau)\|_{L^{2}(\R)}^{1+\frac{1}{s}} \ d\tau
\\[5pt]
& \quad + \frac{4}{\sqrt{\pi}} |\alpha| \ \ve^{b} \|u^{\ve}_{0}\|_{L^{2}(\R)}^{1/2} \int_{0}^{t} \|\partial_x v^{\ve}(\tau)\|_{L^{2}(\R)} \ \|\partial_x u^{\ve}(\tau)\|_{L^{2}(\R)}^{3/2} \ d\tau
  + |\alpha|^{2} \int_{\R} |v^{\ve}(t,x)|^{2} \ dx.
\end{aligned}
\end{equation}

\medskip
3. Now, replacing \eqref{ecuaD014} in \eqref{auxiliar}, we have
$$
\begin{aligned}
&  \int_{\R} |(-\Delta)^{s/2} u^{\ve}(t,x)|^{2} \ dx
+ \ve^{a} \int_{\R} |\partial_x u^{\ve}(t,x)|^{2} \ dx +\frac{1}{4} \int_{\R} |u^{\ve}(t,x)|^{4} \ dx 
\\[5pt]
& \quad  \leq \|(-\Delta)^{s/2} u^{\ve}_{0}\|_{L^{2}(\R)}^{2}
+ \ve^{a} \ \|\partial_x u^{\ve}_{0}\|_{L^{2}(\R)}^{2} + \frac{1}{2} \|u^{\ve}_{0}\|_{L^{4}(\R)}^{4} + \|u^{\ve}_{0}\|_{L^{\infty}(\R)} \|v^{\ve}_{0} \|_{L^{2}(\R)} \|u^{\ve}_{0}\|_{L^{2}(\R)}
 \\[5pt]
& \quad + \frac{4 \ |\alpha|}{\sqrt{\pi(2s-1)}} \ \|g^{\prime}_{\ve}\|_{L^{\infty}(\R)} \ \|u^{\ve}_{0}\|_{L^{2}(\R)}^{1-\frac{1}{2s}}\int_{0}^{t} \| v^{\ve}(\tau)\|_{L^{2}(\R)} \ \|(-\Delta)^{s/2} u^{\ve}(\tau)\|_{L^{2}(\R)}^{1+\frac{1}{2s}} \ d\tau
 \\[5pt]
& \quad + \frac{8 \ |\alpha| \ |\beta|}{\pi(2s-1)} \ \|u^{\ve}_{0}\|_{L^{2}(\R)}^{3-\frac{1}{s}} \int_{0}^{t} \|(-\Delta)^{s/2} u^{\ve}(\tau)\|_{L^{2}(\R)}^{1+\frac{1}{s}} \ d\tau
\\[5pt]
& \quad + \frac{4}{\sqrt{\pi}} |\alpha| \ \ve^{b} \|u^{\ve}_{0}\|_{L^{2}(\R)}^{1/2} \int_{0}^{t} \|\partial_x v^{\ve}(\tau)\|_{L^{2}(\R)} \ \|\partial_x u^{\ve}(\tau)\|_{L^{2}(\R)}^{3/2} \ d\tau,
 \\[5pt]
& \quad + |\alpha|^{2} e^{T} \|v^{\ve}_{0}\|^{2}_{L^{2}(\R)}
  + \frac{16 |\alpha|^{2} \beta^{2} e^{T}}{\pi(2s-1)} \|u^{\ve}_{0}\|_{L^{2}(\R)}^{2-\frac{1}{s}}\int_{0}^{t} \|(-\Delta)^{s/2} u^{\ve}(\tau)\|_{L^{2}(\R)}^{2+\frac{1}{s}} \ d\tau
\end{aligned}
$$
or conveniently we write 
\begin{equation}
\label{ecuaI17}
\begin{aligned}
&  1 + \int_{\R} |(-\Delta)^{s/2} u^{\ve}(t,x)|^{2} \ dx
+ \ve^{a} \int_{\R} |\partial_x u^{\ve}(t,x)|^{2} \ dx + \frac{1}{4} \int_{\R} |u^{\ve}(t,x)|^{4} \ dx \leq \theta(t):= 1
 \\[5pt]
& \quad + \|(-\Delta)^{s/2} u^{\ve}_{0}\|_{L^{2}(\R)}^{2}
+ \ve^{a} \ \|\partial_x u^{\ve}_{0}\|_{L^{2}(\R)}^{2}+\frac{1}{2} \|u^{\ve}_{0}\|_{L^{4}(\R)}^{4} 
\\[5pt]
& \quad + \|u^{\ve}_{0}\|_{L^{\infty}(\R)} 
\|v^{\ve}_{0} \|_{L^{2}(\R)} \|u^{\ve}_{0}\|_{L^{2}(\R)} + |\alpha|^{2} e^{T} \|v^{\ve}_{0}\|^{2}_{L^{2}(\R)}
 \\[5pt]
& \quad + \frac{4 \ |\alpha|}{\sqrt{\pi(2s-1)}} \ \|g^{\prime}_{\ve}\|_{L^{\infty}(\R)} \ \|u^{\ve}_{0}\|_{L^{2}(\R)}^{1-\frac{1}{2s}} 
\int_{0}^{t} \| v^{\ve}(\tau)\|_{L^{2}(\R)} \ \|(-\Delta)^{s/2} u^{\ve}(\tau)\|_{L^{2}(\R)}^{1+\frac{1}{2s}} \ d\tau
 \\[5pt]
& \quad + \frac{8 \ |\alpha| \ |\beta|}{\pi(2s-1)} \ \|u^{\ve}_{0}\|_{L^{2}(\R)}^{3-\frac{1}{s}} \int_{0}^{t} \|(-\Delta)^{s/2} u^{\ve}(\tau)\|_{L^{2}(\R)}^{1+\frac{1}{s}} \ d\tau
\\[5pt]
& \quad + \frac{4}{\sqrt{\pi}} |\alpha| \ \ve^{b} \|u^{\ve}_{0}\|_{L^{2}(\R)}^{1/2} \int_{0}^{t} \|\partial_x v^{\ve}(\tau)\|_{L^{2}(\R)} \ \|\partial_x u^{\ve}(\tau)\|_{L^{2}(\R)}^{3/2} \ d\tau,
 \\[5pt]
& \quad + \frac{16 |\alpha|^{2} \beta^{2} e^{T}}{\pi(2s-1)} \|u^{\ve}_{0}\|_{L^{2}(\R)}^{2-\frac{1}{s}}\int_{0}^{t} \|(-\Delta)^{s/2} u^{\ve}(\tau)\|_{L^{2}(\R)}^{2+\frac{1}{s}} \ d\tau.
\end{aligned}
\end{equation}
From the above definition, we have
$$
\begin{aligned}
 \theta^{\prime}(t)
& \leq \frac{4 \ |\alpha|}{\sqrt{\pi (2s-1)}} \ \|g^{\prime}_{\ve}\|_{L^{\infty}(\R)} \ \|u^{\ve}_{0}\|_{L^{2}(\R)}^{1-\frac{1}{2s}} 
\ \| v^{\ve}(t) \ \|_{L^{2}(\R)} \ \theta(t)^{\frac{1}{2}+\frac{1}{4s}}
\\[5pt]
& + \frac{8 \ |\alpha| \ |\beta|}{\pi(2s-1)} \ \|u^{\ve}_{0}\|_{L^{2}(\R)}^{3-\frac{1}{s}} \ \theta(t)^{\frac{1}{2}+\frac{1}{2s}}
 \\[5pt]
& + \frac{4 |\alpha| \ \ve^{b}}{\sqrt{\pi} \ve^{3a/4}} \ \|u^{\ve}_{0}\|_{L^{2}(\R)}^{1/2} \ \|\partial_x v^{\ve}(t)\|_{L^{2}(\R)} \ \theta(t)^{3/4}
+ \frac{16 |\alpha|^{2} \beta^{2} e^{T}}{\pi(2s-1)} \|u^{\ve}_{0}\|_{L^{2}(\R)}^{2-\frac{1}{s}} \ \theta(t)^{1+\frac{1}{2s}},
\end{aligned}
$$
where we have used \eqref{ecuaI17}. Since $1/2 < s < 1$, then 
$$
   \frac{3}{4} < \frac{1}{2} + \frac{1}{4s} < 1, \quad 1 < \frac{1}{2} + \frac{1}{2s} < \frac{3}{2},  
$$
and consequently dividing the above inequality
by  $\theta(t)^{\frac{1}{2}+\frac{1}{4s}}$, we obtain
$$
\begin{aligned}
& \frac{1}{\frac{1}{2}-\frac{1}{4s}} \ [\theta(t)^{\frac{1}{2}-\frac{1}{4s}}]^{\prime} e^{t}
\leq \frac{4 \ |\alpha|}{\sqrt{\pi (2s-1)}} \ \|g^{\prime}_{\ve}\|_{L^{\infty}(\R)} \ \|u^{\ve}_{0}\|_{L^{2}(\R)}^{1-\frac{1}{2s}} \ \| v^{\ve}(t)\|_{L^{2}(\R)} \ e^{t}
\\[5pt]
& + \frac{8 \ |\alpha| \ |\beta|}{\pi (2s-1)} \ \|u^{\ve}_{0}\|_{L^{2}(\R)}^{3-\frac{1}{s}} \ \theta(t)^{\frac{1}{4s}} \ e^{t}
  \\[5pt]
& + \frac{4}{\sqrt{\pi}} |\alpha| \ \ve^{b} \ \ve^{-3a/4} \ \|u^{\ve}_{0}\|_{L^{2}(\R)}^{1/2} \ \|\partial_x v^{\ve}(t)\|_{L^{2}(\R)} \ e^{t}
+ \frac{16 |\alpha|^{2} \beta^{2} e^{T}}{\pi(2s-1)} \|u^{\ve}_{0}\|_{L^{2}(\R)}^{2-\frac{1}{s}} \ \theta(t)^{\frac{1}{2} + \frac{1}{4s}} \ e^{t},
\end{aligned}
$$
where we have multiplied the inequality by $e^t$. Then, integrating from 0 to $t> 0$
$$
\begin{aligned}
& \int_{0}^{t} \big[ \theta(\tau)^{\frac{1}{2}-\frac{1}{4s}} \big]^{\prime} \ e^{\tau} \ d\tau
\leq \frac{ |\alpha| (4s-2)}{2 s \sqrt{\pi(2s-1)}} \ \|g^{\prime}_{\ve}\|_{L^{\infty}(\R)} \ \|u^{\ve}_{0}\|_{L^{2}(\R)}^{1-\frac{1}{2s}} \int_{0}^{t} \| v^{\ve}(\tau)\|_{L^{2}(\R)} \ e^{\tau} \ d\tau
  \\[5pt]
& + \frac{2 \ |\alpha| \ |\beta|}{s \pi } \ \|u^{\ve}_{0}\|_{L^{2}(\R)}^{3-\frac{1}{s}} \int_{0}^{t} \theta(\tau)^{\frac{1}{4s}} \ e^{\tau} \ d\tau
 + \frac{4 |\alpha|^{2} \beta^{2} e^{T}}{s \pi} \|u^{\ve}_{0}\|_{L^{2}(\R)}^{2-\frac{1}{s}} \ \int_{0}^{t} \theta(\tau)^{\frac{1}{2}+\frac{1}{4s}} \ e^{\tau} \ d\tau
\\[5pt]
& + \frac{|\alpha| \ \ve^{b} \ \ve^{-3a/4}(2s-1)}{\sqrt{\pi} s} \ \|u^{\ve}_{0}\|_{L^{2}(\R)}^{1/2} \int_{0}^{t} \|\partial_x v^{\ve}(\tau)\|_{L^{2}(\R)} \ e^{\tau} \ d\tau
\end{aligned}
$$
and integrating by parts in the left hand side
\begin{equation}
\label{AUX}
\begin{aligned}
&  \theta(t)^{\frac{1}{2}-\frac{1}{4s}} \ e^{t}
\leq \theta(0)^{\frac{1}{2}-\frac{1}{4s}}
+ \int_{0}^{t} \theta(\tau)^{\frac{1}{2}-\frac{1}{4s}} \ e^{\tau} \ d\tau
\\[5pt]
&\quad + \frac{|\alpha| \ (2s-1)}{s \sqrt{2\pi (2s-1)}} \ \|g_\ve^{\prime}\|_{L^{\infty}(\R)} \ \|u^{\ve}_{0}\|_{L^{2}(\R)}^{1-\frac{1}{2s}} \ \big( e^{2t}-1 \big)^{1/2}
 \\[5pt]
& \quad \quad \times \big( e^{T}  \|v^{\ve}_{0}\|^{2}_{L^{2}(\R)}  
+ \frac{16 \beta^{2} e^{T}}{\pi(2s-1)} \|u^{\ve}_{0}\|_{L^{2}(\R)}^{2-\frac{1}{s}}\int_{0}^{t} \|(-\Delta)^{s/2} u^{\ve}(\tau)\|_{L^{2}(\R)}^{2+\frac{1}{s}} \ d\tau \big)^{1/2}
\\[5pt]
& \quad + \frac{2 \ |\alpha| \ |\beta|}{s \pi} \ \|u^{\ve}_{0}\|_{L^{2}(\R)}^{3-\frac{1}{s}} \int_{0}^{t} \theta(\tau)^{\frac{1}{4s}} \ e^{\tau} \ d\tau
\\[5pt]
& \quad + \frac{|\alpha| \ \ve^{b} \ \ve^{-3a/4} \ (2s-1)}{\sqrt{2\pi} \ s} \ \|u^{\ve}_{0}\|_{L^{2}(\R)}^{1/2}  \big( e^{2t}-1 \big)^{1/2}
\\[5pt]
& \quad \quad \times \bigg( \frac{\|v^{\ve}_{0}\|^{2}_{L^{2}(\R)}}{2 \ve^{b}}
+ \frac{8 C_{1,s} \ \beta^{2}}{\ve^{b+1} \pi (2s-1)} \ \|u^{\ve}_{0}\|_{L^{2}(\R)}^{2-\frac{1}{s}}
\int_{0}^{t} \|(-\Delta)^{s/2} u^{\ve}(\tau)\|_{L^{2}(\R)}^{2+\frac{1}{s}} \ d\tau \bigg)^{1/2}
\\[5pt]
& \quad + \frac{4 |\alpha|^{2} \beta^{2} e^{T} }{s\pi} \|u^{\ve}_{0}\|_{L^{2}(\R)}^{2-\frac{1}{s}} \ \int_{0}^{t} \theta(\tau)^{\frac{1}{2}+\frac{1}{4s}} \ e^{\tau} \ d\tau,
\end{aligned}
\end{equation}
where we have used Holder's inequality and equations \eqref{ecuaI15}-\eqref{ecuaD014}. 

\medskip
4. The goal now is to apply the Generalized Gronwall Lemma (Section \ref{GGL}).
We observe that
$$
\begin{aligned}
   \theta(0)= \bigg[ 1
&+ \|(-\Delta)^{s/2} u^{\ve}_{0}\|_{L^{2}(\R)}^{2}
+ \ve^{a} \ \|\partial_x u^{\ve}_{0}\|_{L^{2}(\R)}^{2}
+ \frac{1}{2} \|u^{\ve}_{0}\|_{L^{4}(\R)}^{4}
\\[5pt]
& + \|u^{\ve}_{0}\|_{L^{\infty}(\R)} \|v^{\ve}_{0} \|_{L^{2}(\R)} \|u^{\ve}_{0}\|_{L^{2}(\R)} 
+ |\alpha|^{2} e^{T} \|v^{\ve}_{0}\|^{2}_{L^{2}(\R)} \bigg],
\end{aligned}
$$
hence from that and taking the square in equation \eqref{AUX}, we have
$$
\begin{aligned}
&  \theta(t)^{1-\frac{1}{2s}} \ e^{2 t}
\leq 2^6 \bigg[ 1
+ \|(-\Delta)^{s/2} u^{\ve}_{0}\|_{L^{2}(\R)}^{2}
+ \ve^{a} \ \|\partial_x u^{\ve}_{0}\|_{L^{2}(\R)}^{2}
+ \frac{1}{2} \|u^{\ve}_{0}\|_{L^{4}(\R)}^{4}
\\[5pt]
& \quad + \|u^{\ve}_{0}\|_{L^{\infty}(\R)} \|v^{\ve}_{0} \|_{L^{2}(\R)} \|u^{\ve}_{0}\|_{L^{2}(\R)} 
+ |\alpha|^{2} e^{T} \|v^{\ve}_{0}\|^{2}_{L^{2}(\R)} \bigg]^{1-\frac{1}{2s}} 
\\[5pt]
& \quad + 2^6 t^{2} \bigg( \fint_0^t \theta(\tau)^{\frac{1}{2}-\frac{1}{4 s}} \ e^{\tau} \ d\tau \bigg)^{2}
+ \frac{ 2^5 |\alpha|^{2} \ (2s-1)}{s^{2} \pi} \ \|g_\ve^{\prime}\|_{L^{\infty}(\R)}^{2} \ \|u^{\ve}_{0}\|_{L^{2}(\R)}^{2-\frac{1}{s}} \ e^{2 t} 
\\[5pt]
& \quad \quad \times \bigg( e^{T}  \|v^{\ve}_{0}\|^{2}_{L^{2}(\R)}  
+ \frac{16 \beta^{2} e^{T}}{\pi(2s-1)} \|u^{\ve}_{0}\|_{L^{2}(\R)}^{2-\frac{1}{s}}\int_{0}^{t} \|(-\Delta)^{s/2} u^{\ve}(\tau)\|_{L^{2}(\R)}^{2+\frac{1}{s}} \ d\tau \bigg)
\\[5pt]
& \quad + \frac{2^8 |\alpha|^{2} \ |\beta|^{2} t^{2}}{s^{2} \pi^{2} } \ \|u^{\ve}_{0}\|_{L^{2}(\R)}^{6-\frac{2}{s}} \bigg(  \fint_0^t \theta(\tau)^{\frac{1}{4 s}} \ e^{\tau} \ d\tau \bigg)^{2}
  \\[5pt]
& \quad  + \frac{2^5 |\alpha|^{2} \ \ve^{2b} \ \ve^{-3a/2} \ (2s-1)^{2}}{\pi s^{2}} \ \|u^{\ve}_{0}\|_{L^{2}(\R)} \ e^{2 t} 
\\[5pt]
&\quad \quad  \times \bigg( \frac{\|v^{\ve}_{0}\|^{2}_{L^{2}(\R)}}{2 \ve^{b}}
+ \frac{8 C_{1,s} \ |\beta|^{2}}{\ve^{b+1} \pi (2s-1)} \ \|u^{\ve}_{0}\|_{L^{2}(\R)}^{2-\frac{1}{s}}
\int_{0}^{t} \|(-\Delta)^{s/2} u^{\ve}(\tau)\|_{L^{2}(\R)}^{2+\frac{1}{s}} \ d\tau \bigg)
  \\[5pt]
& \quad + \frac{2^{10} |\alpha|^{4} |\beta|^{4} e^{2T} }{s^{2}\pi^{2}} \|u^{\ve}_{0}\|_{L^{2}(\R)}^{4-\frac{2}{s}} \ t^2 
\bigg( \fint_0^t \theta(\tau)^{\frac{1}{2}+\frac{1}{4 s}} \ e^{\tau} \ d\tau \bigg)^2.
\end{aligned}
$$
Then, we apply Jesen's inequality to obtain  
$$
\begin{aligned}
&  \theta(t)^{1-\frac{1}{2s}} \ e^{2 t}
\leq 2^6 \bigg[ 1
+ \|(-\Delta)^{s/2} u^{\ve}_{0}\|_{L^{2}(\R)}^{2}
+ \ve^{a} \ \|\partial_x u^{\ve}_{0}\|_{L^{2}(\R)}^{2}
+ \frac{1}{2} \|u^{\ve}_{0}\|_{L^{4}(\R)}^{4}
\\[5pt]
& \quad + \|u^{\ve}_{0}\|_{L^{\infty}(\R)} \|v^{\ve}_{0} \|_{L^{2}(\R)} \|u^{\ve}_{0}\|_{L^{2}(\R)} 
+ |\alpha|^{2} e^{T} \|v^{\ve}_{0}\|^{2}_{L^{2}(\R)} \bigg]^{1-\frac{1}{2s}} 
\\[5pt]
& \quad + 2^6 T \int_0^t \theta(\tau)^{1-\frac{1}{2 s}} \ e^{ 2 \tau} \ d\tau
+ \frac{ 2^5 |\alpha|^{2} \ (2s-1)}{s^{2} \pi} \ \|g_\ve^{\prime}\|_{L^{\infty}(\R)}^{2} \ \|u^{\ve}_{0}\|_{L^{2}(\R)}^{2-\frac{1}{s}} \ e^{2 T} 
\\[5pt]
& \quad \quad \times \bigg( e^{T}  \|v^{\ve}_{0}\|^{2}_{L^{2}(\R)}  
+ \frac{16 \beta^{2} e^{T}}{\pi(2s-1)} \|u^{\ve}_{0}\|_{L^{2}(\R)}^{2-\frac{1}{s}}\int_{0}^{t} \|(-\Delta)^{s/2} u^{\ve}(\tau)\|_{L^{2}(\R)}^{2+\frac{1}{s}} \ d\tau \bigg)
\\[5pt]
& \quad + \frac{2^8 |\alpha|^{2} \ |\beta|^{2} T}{s^{2} \pi^{2} } \ \|u^{\ve}_{0}\|_{L^{2}(\R)}^{6-\frac{2}{s}} \int_0^t \theta(\tau)^{\frac{1}{2 s}} \ e^{2 \tau} \ d\tau
  \\[5pt]
& \quad  + \frac{2^5 |\alpha|^{2} \ \ve^{2b} \ \ve^{-3a/2} \ (2s-1)^{2}}{\pi s^{2}} \ \|u^{\ve}_{0}\|_{L^{2}(\R)} \ e^{2 T} 
\\[5pt]
&\quad \quad  \times \bigg( \frac{\|v^{\ve}_{0}\|^{2}_{L^{2}(\R)}}{2 \ve^{b}}
+ \frac{8 C_{1,s} \ |\beta|^{2}}{\ve^{b+1} \pi (2s-1)} \ \|u^{\ve}_{0}\|_{L^{2}(\R)}^{2-\frac{1}{s}}
\int_{0}^{t} \|(-\Delta)^{s/2} u^{\ve}(\tau)\|_{L^{2}(\R)}^{2+\frac{1}{s}} \ d\tau \bigg)
  \\[5pt]
& \quad + \frac{2^{10} |\alpha|^{4} |\beta|^{4} e^{2T} }{s^{2}\pi^{2}} \|u^{\ve}_{0}\|_{L^{2}(\R)}^{4-\frac{2}{s}} 
\ T \int_0^t \theta(\tau)^{1+\frac{1}{2 s}} \ e^{2 \tau} \ d\tau.
\end{aligned}
$$
Moreover, after an algebraic manipulation and using that $e^t> 1$ for any $t>0$, we may write
\begin{equation}
\label{INIGW}
\begin{aligned}
 \theta(t)^{1-\frac{1}{2s}} \ e^{2 t}
\leq C &+ C_1 \int_0^t \theta(\tau)^{1-\frac{1}{2 s}} \ e^{ 2 \tau} \ d\tau
\\[5pt]
& + C_2 \int_0^t \theta(\tau)^{\frac{1}{2 s}} \ e^{2 \tau} \ d\tau + C_3 \int_{0}^{t} \theta(\tau)^{1+\frac{1}{2 s}} e^{2 \tau}\ d\tau, 
\end{aligned}
\end{equation}
where 
$$
\begin{aligned} 
 &  C:= 2^6 \big( 1
+ \|(-\Delta)^{s/2} u^{\ve}_{0}\|_{L^{2}(\R)}^{2}
+ \|\partial_x u^{\ve}_{0}\|_{L^{2}(\R)}^{2}
+ \frac{1}{2} \|u^{\ve}_{0}\|_{L^{4}(\R)}^{4}
\\[5pt]
& \quad + \|u^{\ve}_{0}\|_{L^{\infty}(\R)} \|v^{\ve}_{0} \|_{L^{2}(\R)} \|u^{\ve}_{0}\|_{L^{2}(\R)} 
+ |\alpha|^{2} e^{T} \|v^{\ve}_{0}\|^{2}_{L^{2}(\R)} \big)^{1-\frac{1}{2s}} 
\\[5pt]
& \quad + \frac{ 2^5 |\alpha|^{2} \ (2s-1)}{s^{2} \pi} \ \|g_\ve^{\prime}\|_{L^{\infty}(\R)}^{2} \ \|u^{\ve}_{0}\|_{L^{2}(\R)}^{2-\frac{1}{s}}  \|v^{\ve}_{0}\|_{L^{2}(\R)}^{2} \ e^{3 T}
\\[5pt]
& \quad + \frac{2^4 |\alpha|^{2} \ \ve^{b} \ \ve^{-3a/2} \ (2s-1)^{2}}{\pi s^{2}} \ \|u^{\ve}_{0}\|_{L^{2}(\R)} \|v^{\ve}_{0}\|^{2}_{L^{2}(\R)} \ e^{2 T},
\\[5pt]
& C_1:= 2^6 T,
\quad \quad  C_2:=  \frac{2^8 |\alpha|^{2} \ |\beta|^{2} T}{s^{2} \pi^{2} } \ \|u^{\ve}_{0}\|_{L^{2}(\R)}^{6-\frac{2}{s}}, 
\\[5pt]
&C_3:= \frac{ 2^9 |\alpha|^{2} \ |\beta|^{2}}{ s^{2} \pi^2} \ \|g_\ve^{\prime}\|_{L^{\infty}(\R)}^{2} \ \|u^{\ve}_{0}\|_{L^{2}(\R)}^{4-\frac{2}{s}} \ e^{3 T}
  \\[5pt]
& \quad  + \frac{2^8 C_{1,s} |\alpha|^{2} \ |\beta|^{2} \ \ve^{b-1} \ \ve^{-3a/2} \ (2s-1)}{ \pi^2 s^{2}} \ \|u^{\ve}_{0}\|_{L^{2}(\R)}^{3-\frac{1}{s}} \ e^{2 T} 
+ \frac{2^{10} |\alpha|^{4} \ |\beta|^{4} e^{2T}}{\pi^2 s^2} \|u^{\ve}_{0}\|_{L^{2}(\R)}^{4-\frac{2}{s}} \ T.
\end{aligned}
$$
Therefore, taking $a= 4$ and $b= 7$ the above positive constants $C$, $C_1$, $C_2$ and $C_3$ are independent of $\ve> 0$.  
Now, since 
$$
    (1-\frac{1}{2 s})(\frac{2s+1}{2s-1}) = 1+\frac{1}{2 s}, \quad \text{and} \quad (1-\frac{1}{2 s})(\frac{1}{2s-1}) = \frac{1}{2 s},
$$ 
then we have from \eqref{INIGW} 
$$
\begin{aligned}
 \theta(t)^{1-\frac{1}{2s}} \ e^{2 t}
\leq C &+ C_1 \int_0^t \theta(\tau)^{1-\frac{1}{2 s}} \ e^{ 2 \tau} \ d\tau
\\[5pt]
& + C_2 \int_0^t \big(\theta(\tau)^{1-\frac{1}{2 s}}\big)^{\frac{1}{2s-1}} \ e^{2 \tau} \ d\tau 
+ C_3 \int_{0}^{t} \big(\theta(\tau)^{1-\frac{1}{2 s}}\big)^{\frac{2s+1}{2s-1}} e^{2 \tau}\ d\tau.
\end{aligned}
$$
For each $1/2 < s < 1$ we have 
$$
    1 < \frac{1}{2s-1} < \frac{2s+1}{2s-1}
$$ 
therefore from the above inequality we may write
$$
\begin{aligned}
 \theta(t)^{1-\frac{1}{2s}} \ e^{2 t}
\leq C &+ C_1 \int_0^t \theta(\tau)^{1-\frac{1}{2 s}} \ e^{ 2 \tau} \ d\tau
\\[5pt]
& + C_2 \int_0^t \big(\theta(\tau)^{1-\frac{1}{2 s}} \ e^{2 \tau} \big)^{\frac{2s+1}{2s-1}} \ d\tau 
+ C_3 \int_{0}^{t} \big(\theta(\tau)^{1-\frac{1}{2 s}} e^{2 \tau} \big)^{\frac{2s+1}{2s-1}}  \ d\tau
\end{aligned}
$$
or defining $\eta(t):= \theta(t)^{1-\frac{1}{2 s}} \ e^{2 t}$
\begin{equation}
 \eta(t) \leq C + \int_0^t \bigg[ C_1 \eta(\tau) + \Big(C_2 + C_3\Big) \big(\eta(\tau) \big)^{\frac{2s+1}{2s-1}} \bigg] \ d\tau.
\end{equation}
Therefore, applying the Generalized Grownwall Lemma, more precisely \eqref{ecuaI22} with $\sigma= \frac{2s+1}{2s-1}> 1$, 
we must have for each $s \in (1/2,1)$ 
$$
\begin{aligned}
  C 
  & < \Big\{ \exp \bigg[ \Big( 1 - \frac{2s+1}{2s-1} \Big) \int_{0}^{T} C_1 \ d\tau \big] \Big\}^{\frac{1}{\frac{2s+1}{2s-1}-1}} \
   \Big\{\Big( \frac{2s+1}{2s-1} - 1 \Big) \int_{0}^{T} \Big(C_2 + C_3\Big) \ d\tau \Big\}^{-\frac{1}{\frac{2s+1}{2s-1}-1}}
    \\[5pt]
   & \quad = \frac{(2s-1)^{\frac{2s-1}{2}} \ \exp \big[ - C_1 T \big]}{\Big \{2 \Big(C_2 + C_3\Big) \ T \Big\}^{\frac{2s-1}{2}}}
\end{aligned}
$$
or equivalently 
\begin{equation}
\label{conditionC}
  C \ \Big(C_2 + C_3\Big)^{\frac{2s-1}{2}} \exp[64 T^2] \ T^{\frac{2s-1}{2}} \leq \Big({\frac{2s-1}{2}}\Big)^{\frac{2s-1}{2}}.
\end{equation}
One remarks that 
$$
    \lim_{s\rightarrow\frac{1}{2}} \Big(\frac{2s-1}{2}\Big)^{\frac{2s-1}{2}}= 1.
$$
Hence for any $s \in (1/2,1)$ fixed, there exists $\alpha_0> 0$ and $E_0> 0$, such that 
condition \eqref{conditionC} is satisfied when $\|u_{0}\|_{L^{2}(\R)} \leq E_0$, or $|\alpha| \leq \alpha_0$. In fact, if there is no coupling, that is
$\alpha= 0$ ($C_2= C_3= 0$), then condition \eqref{conditionC} is trivially satisfied.  
Consequently, we have
$$
\begin{aligned}
  \eta(t)
   & \leq C \ \bigg\{ \exp \bigg[ \Big( 1-\frac{2s+1}{2s-1} \Big) \int_{0}^{t} C_1 \ d\tau \bigg]
 \\[5pt]
& - C^{-1} \ \Big( \frac{2s+1}{2s-1}-1 \Big) \int_{0}^{t} \Big(C_2 + C_3\Big) \ \exp \bigg[ \Big( 1 - \frac{2s+1}{2s-1} \Big) \int_{\tau}^{t} C_1 \ dr \bigg] \ d\tau \bigg\}^{\frac{1}{\frac{2s+1}{2s-1}-1}}
 \\[5pt]
   & = C \ \bigg\{ \exp \bigg[ \frac{2 C_1}{1-2s} t \big]
   - C^{-1} \frac{2}{2s-1} \ \Big(C_2 + C_3\Big) \int_{0}^{t} \exp \big[ \frac{2 C_1}{1-2s} (t - \tau) \big] \ d\tau \bigg\}^{\frac{2s-1}{2}}
  \\[5pt]
& = C \ \bigg\{ \exp \bigg[ \frac{2 C_1}{1-2s} t \bigg] - \frac{C^{-1} \Big(C_2 + C_3\Big)}{C_1}\bigg(1-\exp \bigg[ \frac{2 C_1}{1-2s} t \bigg]\bigg) \bigg\}^{\frac{2s-1}{2}},
\end{aligned}
$$
from which follows the proof of the theorem.
\end{proof}

\medskip
Finally, we establish a maximum principle for the solution $v^\ve$ of $(\ref{ecuaI2})_2$. 
\begin{proposition}[Maximum Principle]
\label{MaxPrincip}
Let $(u^\ve,v^\ve)$ be the unique solution of \eqref{ecuaI2}. Then, $v^\ve$ satisfies
\begin{equation}
\label{MaxPrincip20}
   \sup_{(0,T) \times \R} |v^\ve| \leq  \|v_0\|_{L^\infty(\R)}.
\end{equation}
\end{proposition}
\begin{proof}
For $\ve> 0$ fixed, let us define $w:= v^\ve - \|v_0\|_{L^\infty(\R)}$, and we 
will show that, $w^+= \max \{w,0\}= 0$. Clearly $w^+ \geq 0$, then we assume by 
contradiction that, $w^+ > 0$. Therefore, there exists $\mu> 0$, such that $w^+ \ge \mu$. 
Since $w^+(t) \in H^1(\R)$ for each $t \in (0,T)$, we can use $w^+$ as a test function for equation $(\ref{ecuaI2})_2$, 
(similar to equation \eqref{RPA2}), that is to say
$$
\begin{aligned}
\frac{d}{dt}\!\! \int_\R |w^+(t)|^2 \, dx &+ \!\! \int_\R g_\ve(v^\ve(t)) \, (-\Delta)^{s/2} w^+(t) \ dx
\\[5pt]
&= \beta \!\!  \int_\R |u^\ve(t)|^2 \, (-\Delta)^{s/2} w^+(t) \ dx
-  \!\! \int_\R \ve^7 \,  |\nabla w^+(t)|^2 \ dx. 
\end{aligned}
$$
Now, we consider the following estimate 
$$
\begin{aligned}
g_\ve(v^\ve) \, (-\Delta)^{s/2} w^+ &+ \ve^7 \, |\nabla w^+|^2 -  \beta \,  |u^\ve|^2 \, (-\Delta)^{s/2} w^+
\\[5pt]
&\geq - M \, |w^+|  \, |(-\Delta)^{1/2} w^+| + \ve^7 \, |\nabla w^+|^2 - \frac{|\beta| \, \|u^\ve\|^2_{L^\infty(\R)}}{\mu} |(-\Delta)^{1/2} w^+|
\\[5pt]
&\geq  -\frac{ \big(M \, \mu + |\beta| \, \|u^\ve\|^2_{L^\infty(\R)}\big)^2}{2 \, \ve^7 \mu^2} |w^+|^2 - \ve^7 \,  |(-\Delta)^{1/2} w^+|^2 + \ve^7 \, |\nabla w^+|^2,
\end{aligned}
$$
where we have used Young's inequality. 
Consequently, since $w^+(0)= 0$ we have 
$$
    \int_\R |w^+(t)|^2 \, dx \leq  \frac{ \big(M \, \mu + |\beta| \, \|u^\ve\|^2_{L^\infty(\R)}\big)^2}{2 \, \ve^7 \mu^2} \int_0^t \!\! \! \int_\R |w^+(\tau,x)|^2 \, dx \, d\tau, 
$$
which applying the Gronwall's lemma implies a contradiction. Therefore, we have 
\begin{equation}
\label{MaxPrincip10}
   w^+(t)=(v^\ve(t) - \|v_0\|_{L^\infty(\R)})^+ \equiv 0.
\end{equation}   
A similar argument can 
also show that
\begin{equation}
\label{MaxPrincip15}
   (-v^\ve(t) - \|v_0\|_{L^\infty(\R)})^+ \equiv 0.
\end{equation}   
Nonetheless, the equations \eqref{MaxPrincip10} and \eqref{MaxPrincip15} mean nothing other than   
$$
    \|v^\ve(t)\|_{L^\infty(\R)} \leq \|v_0\|_{L^\infty(\R)}  \quad \text{for each $t \in (0,T)$}. 
$$
\end{proof}

\section{Existence of Weak Solutions}

The main issue of this section is to show 
the solvability of the Cauchy problem \eqref{ecua00},
that is, we prove Theorem \ref{MAINTHM} (Main Theorem). More 
precisely, from the equivalence of mild solutions (when it exists) and weak 
solutions, we obtain a weak formulation from \eqref{ecua152},
see Lemma \ref{EquivMildWeak}, and the goal is to pass to the limit 
as $\ve \to 0^+$ to show a solution of the 
Cauchy problem \eqref{ecua00} in the sense of Definition \ref{DEFSOL}.
We apply the Aubin-Lions Theorem 
to show that the family $\{u^\ve\}$ is relatively compact in $L^2$. 
The similar result for the family $\{v^\ve\}$ does not follow 
analogously, since $(\ref{ecua00})_2$ degenerates. Hence
we apply Tartar's methodology in \cite{Tartar2}, (see also \cite{Malek}),
adapted to our context of fractional porous media equation.  

\medskip
First, we have the following 
\begin{lemma}
\label{EquivMildWeak}
Let $\alpha_0> 0, E_0> 0$ be given by Theorem \ref{Second estimate}, such that, $|\alpha| \leq \alpha_0$ or
$\|u_{0}\|_{L^{2}(\R)} \leq E_0$. Then, 
the unique mild solution $(u^{\ve},v^{\ve})$ of \eqref{ecuaI2} satisfies, 
\begin{equation}
\label{RPA1}
\begin{aligned}
    & i \!\!\int_{0}^{T}\int_{\R} \Big( u^{\ve}(t,x) \ \partial_t \overline{\varphi}(t,x) 
     + (-\Delta)^{s/2}u^{\ve}(t,x) \ (-\Delta)^{s/2} \overline{\varphi}(t,x) \Big) dx dt
    + i \!\! \int_{\R} u_{0}^{\ve}(x) \ \overline{\varphi}(0,x) dx
\\[5pt]
    &- \ve^{a} \int_{0}^{T} \int_{\R} u^{\ve}(t,x) \ \Delta \overline{\varphi} \ dx dt + \alpha  \int_{0}^{T} \int_{\R} v^{\ve}(t,x) \ u^{\ve}(t,x) \ \overline{\varphi}(t,x) \ dxdt
\\[5pt]
 & \quad + \int_{0}^{T} \int_{\R} |u^{\ve}(t,x)|^{2} \ u^{\ve}(t,x) \ \overline{\varphi}(t,x) \ dx dt = 0,
\end{aligned}
\end{equation}
\begin{equation}
\label{RPA2}
\begin{aligned}
   & \int_{0}^{T} \!\!\! \int_{\R} v^{\ve}(t,x) \ \partial_t \psi(t,x)
   - g_\ve(v^{\ve}(t,x)) \ (-\Delta)^{s/2} \psi(t,x) \ dx dt 
   + \int_{\R} v_{0}^{\ve}(x) \ \psi(0,x) \ dx
   \\[5pt]
    &+ \ve^{7} \int_{0}^{T}\!\!\! \int_{\R} v^{\ve}(t,x) \ \Delta\psi(t,x) \ dx dt + \beta \int_{0}^{T} \!\!\! \int_{\R} |u^{\ve}|^{2}(t,x) \ (-\Delta)^{s/2} \psi(t,x) \ dt dx = 0 
\end{aligned}
\end{equation}
for each test functions $\varphi, \psi \in C^{\infty}_{c} \big( (-\infty,T) \times \R \big)$, with $\varphi$ being complex-valued and $\psi$ real-valued.

Moreover, there exists a positive constant $C$ independent of $\ve> 0$, such that 
\begin{equation}
\label{EQSMHUM}
\begin{aligned}
  \int_0^T \|\partial_t u^{\ve}(t)\|^2_{H^{-1}(\R)} dt \leq C, 
  \qquad
 \int_0^T \|\partial_t v^{\ve}(t)\|^2_{H^{-1}(\R)} dt \leq C.
\end{aligned}
\end{equation}
\end{lemma}

\begin{proof}
Equations \eqref{RPA1}, \eqref{RPA2} are obtained from \eqref{ecua152}, that is, applying
the equivalence between mild solutions and weak solutions, (see Ball \cite{Ball}, p. 371), 
which are obtained via functional analysis arguments. Similarly, 
the inequalities in equation \eqref{EQSMHUM} are obtained from the weak formulation, i.e.
equations \eqref{RPA1} and \eqref{RPA2}, 
applying standard functional analysis results, the uniform
boundedness of $u_{0}^{\ve}$, $v_{0}^{\ve}$, and also the uniform estimates 
from Lemma \ref{First estimate} and Theorem \ref{Second estimate}.  
\end{proof}

\subsection{Proof of main theorem}

Now, we are ready to show the main result of this article. 
\begin{proof}[\bf Proof]
1. First, under the conditions of Lemma \ref{EquivMildWeak}, for each $\ve> 0$, 
let $(u^{\ve},v^{\ve}) \in C([0,T);H^1(\R)) \times C([0,T);H^1(\R))$ be
the unique mild solution of \eqref{ecuaI2}, 
satisfying \eqref{ecua152}. Then, the pair$(u^{\ve}(t,x), v^{\ve}(t,x))$
satisfies the equations \eqref{RPA1} and \eqref{RPA2}.
To obtain \eqref{RP1}, \eqref{RP2} we pass to the limit respectively in 
\eqref{RPA1} and \eqref{RPA2} as $\ve \to 0^{+}$. 
Therefore, we need to show strong convergence, which implies 
a.e. convergence (along subsequences) of the sequences  
$\{u^{\ve}\}_{\ve> 0}$, and $\{v^{\ve}\}_{\ve> 0}$.

\medskip
2. Let us show that the family $\{u^{\ve}\}_{\ve> 0}$ 
is relatively compact. From \eqref{ecuaI9}, \eqref{ecuaI13}, 
it follows that $\{u^{\ve}\}_{\ve> 0}$ is (uniformly) bounded in 
$L^{\infty}(0,T;H^s(\R))$, hence it is possible to select a subsequence,
still denoted by $\{u^{\ve}\}_{\ve> 0}$, which converges weakly-$\star$ to $u$ in $L^{\infty}(0,T;H^s(\R))$. 
Applying the Rellich's Theorem, for any compact set $K \subset \R$,
the embedding of $H^{s}(K)$ in $L^{2}(K)$ is compact. Therefore, 
since the sequence $\{u^{\ve}\}_{\ve> 0}$ is uniformly bounded in 
$L^{2}(0,T;H^{-1}(\R))$, we apply the Aubin-Lions Theorem and 
obtain (along a suitable subsequence) that $u^{\ve}$ converges 
strongly to $u$ in $L^{2}(0,T;L^{2}(K))$, and thus
\begin{equation}
\label{uconvergence}
 \text{$u^{\ve}(t,x) \to u(t,x)$ as $\ve \to 0$ almost everywhere
in $(0,T) \times \R$.}
\end{equation} 

\medskip
3. Now, we show that the family $\{v^{\ve}\}_{\ve> 0}$ is relatively compact. 
First, we multiply equation $(\ref{ecuaI2})_2$ by $\eta_k^\prime(v^\ve)$, (see Section \ref{Entropies}),
 and applying a standard procedure (e.g. the theory of scalar conservation laws), 
 we obtain in distribution sense
\begin{equation}
\label{entropyv}
\begin{aligned}
\partial_t \eta_k(v^\ve) + (-\Delta)^{s/2} |g_\ve(v^\ve) - g_\ve(k)|&
=  \beta \, \eta_k^\prime(v^\ve) (-\Delta)^{s/2} |u^\ve|^2 + \ve^7 \Delta \eta_k(v^\ve)
\\[5pt]
& - \ve^7 |\nabla v^\ve|^2  \eta_k^{\prime \prime}(v^\ve) -  R^\ve_k, 
\end{aligned} 
\end{equation} 
where we have used \eqref{remainder} and obvious notation.
Let $\eta$ be any smooth (say $C^2$) entropy 
(which is linear at infinity, i.e. $\eta^{\prime \prime} (\cdot) \in C_0(\R)$). Then, multiplying equation \eqref{entropyv} by 
$\eta^{\prime \prime}(k)$ and integrating in $\R$ with respect to $k$,
we obtain in the sense of distributions 
\begin{equation}
\label{SENTROPIA}
\begin{aligned}
\partial_t \eta(v^\ve) + (-\Delta)^{s/2} q(v^\ve)& 
= \beta \, \eta^\prime(v^\ve) (-\Delta)^{s/2} |u^\ve|^2 + \ve^7 \Delta \eta(v^\ve)
\\[5pt]
& - \ve \, (-\Delta)^{s/2} \eta(v^\ve) - \ve^7 |\nabla v^\ve|^2  \eta^{\prime \prime}(v^\ve) -  \mathcal{R}^\ve, 
\end{aligned} 
\end{equation}
where the function $q$ satisfies $q^\prime= \eta^\prime g^\prime$, (recall that $g_\ve(v^\ve)= g(v^\ve) + \ve v^\ve)$, and 
$$
\mathcal{R}^\ve:= \frac{1}{2} \int_{\R} \eta^{\prime \prime}(k) \, R^\ve_k \, dk. 
$$
From \eqref{ecuaD014}, \eqref{MaxPrincip20}
it follows that the family $\{v^{\ve}\}_{\ve> 0}$ is (uniformly) bounded in $L^{\infty}(0,T;L^{2}(\R) \cap L^{\infty}(\R))$,
hence it is possible to select a subsequence, still denoted by $\{v^{\ve}\}_{\ve> 0}$, 
which converges weakly-$\star$ to $v$ in $L^{\infty}(0,T;L^{2}(\R) \cap L^{\infty}(\R))$. Moreover, we show that for any entropy pair $(\eta,q)$, 
\begin{equation}
\label{compactnessv}
  \partial_t \eta(v^\ve) + (-\Delta)^{s/2} q(v^\ve) \in \Big\{\text{compact set of $H^{-1}_{\loc}((0,\infty) \times \R)$}\Big\}. 
\end{equation}
Indeed, we first observe that the left hand side of \eqref{SENTROPIA} is uniformly bounded 
in $W^{-1,\infty}_\loc((0,\infty) \times \R)$. From equation \eqref{ecuaI15} the terms 
$\ve^7 \Delta \eta(v^\ve)$,  $\ve \, (-\Delta)^{s/2} \eta(v^\ve)$
are compact in $H^{-1}_{\loc}((0,\infty) \times \R)$,
let us show the former the second one is similar. 
Let $K \subset (0,T) \times \R$ be any compact set, and
$\phi \in C^\infty_c((0,T) \times \R)$. Then, we have
$$
\begin{aligned}
|\langle \ve^{7} \Delta \eta(v^{\ve}), \phi \rangle| &
\leq \int_0^T\!\!\! \int_\R| \ve^{7} \nabla \eta(v^{\ve}) \cdot  \nabla \phi | \, dx dt
\\[5pt]
& \leq \ve^{7/2} \bigg(\int_0^T\!\!\! \int_\R | \ve^{7/2} \nabla \eta(v^{\ve})|^{2} \ dx \ dt \bigg)^{1/2} \| \nabla \phi\|_{L^2(K)}
\leq \ve^{7/2} C \,  \| \nabla \phi\|_{L^2(K)}, 
\end{aligned}
$$
where $C> 0$ does not depend on $\ve$, 
and we have used that $\eta$ is linear at infinity. 
 Then, taking the supremum with respect to the
set $W= \{ \phi \in H^1: \| \phi \|_{H^1} \leq 1\}$ and
passing to the limit as $\ve \to 0^+$, the family 
$\{ \ve^7 \Delta \eta(v^\ve) \}$ converges to 
zero in $H^{-1}_{\loc}$.  
 
From equation equation \eqref{ecuaI15}
the family $\{\ve^7 |\nabla v^\ve|^2  \eta^{\prime \prime}(v^\ve)\}$ is uniformly 
bounded in $L^1_\loc((0,\infty) \times \R)$, and thus in the space of Radon measures $\clg{M}_{\loc}((0,\infty) \times \R)$. 
 Hence compact in $W^{-1,q}_\loc((0,\infty) \times \R)$, for $1 \leq q< 3/2$. Similarly,
 from \eqref{ecuaI13}, \eqref{new}, the family 
 $\{\eta^\prime(v^\ve) (-\Delta)^{s/2} |u^\ve|^2\}$ is uniformly 
bounded in $L^1_\loc((0,\infty) \times \R)$, and hence compact in 
$W^{-1,q}_\loc((0,\infty) \times \R)$, for $1 \leq q< 3/2$. Finally, we consider the 
following 

\smallskip
 \underline {Claim:} The family $\{\clg{R}^\ve \}$ is compact in 
$W^{-1,q}_\loc((0,\infty) \times \R)$, for $1 \leq q< 3/2$. 
 
 \medskip
Consequently, we have  $\partial_t \eta(v^\ve) + (-\Delta)^{s/2} q(v^\ve)$ in 
$$
      \Big\{\text{bounded set of $W^{-1,\infty}_{\loc}((0,\infty) \times \R)$}\Big\}
     \cap \Big\{\text{compact set of $W^{-1,q}_{\loc}((0,\infty) \times \R)$}\Big\},  
$$
and due to a well known interpolation argument, (see Lemma 3.12 in \cite{Malek}), it follows \eqref{compactnessv}. 

\medskip
Proof of Claim: It is enough to show that, the family $\{\clg{R}^\ve \}$ is uniformly bounded in 
$\clg{M}_{\loc}((0,\infty) \times \R)$. To this end, we observe that the left hand side of \eqref{SENTROPIA} is also uniformly bounded 
in $H^{-1}_\loc((0,\infty) \times \R)$, and jointly with the others terms clearly shows that, $\clg{R}^\ve$ is a uniformly 
bounded distribution ( in $\clg{D}^\prime$), which is positive by definition, hence from a well known result a Radon measure. 

\medskip
Now, from \eqref{compactnessv} we may apply the Tartar's method in \cite{Tartar2}, (see also \cite{Malek}), 
which implies the compactness of the sequence $\{v^\ve\}$ in $L^1_\loc((0,\infty) \times \R)$. 
Thus along a suitable subsequence
\begin{equation}
\label{vconvergence}
\text{$v^\ve(t,x) \to v(t,x)$ as $\ve \to 0$ 
almost everywhere in $(0,T) \times \R$.}
\end{equation}

\medskip
4. Finally, from \eqref{uconvergence}, \eqref{vconvergence} and 
due to a standard diagonalization procedure, we apply the Dominated Convergence Theorem 
to pass to the limit as $\ve \to 0$ 
in the equations 
\eqref{RPA1} and \eqref{RPA2}, which togheter the
Definition \ref{DEFSOL}
gives the solvability of the Cauchy problem \eqref{ecua00}. Moreover, inequality \eqref{maxprinciple} 
follows from \eqref{MaxPrincip20}. 
\end{proof}

\section*{Acknowledgements}

Conflict of Interest: Author Wladimir Neves has received research grants from CNPq
through the grant  308064/2019-4, and also by FAPERJ 
(Cientista do Nosso Estado) through the grant E-26/201.139/2021. 


\end{document}